\documentclass{gtart_h}

\def\ifplaintex{\expandafter\ifx\csname documentclass\endcsname\relax}

\def\gtp{{\mathsurround=0pt\it $\cal G\mskip-2mu$eometry \&\ 
$\cal T\!\!$opology $\cal P\!$ublications}}  

\def\recd{{\small Received:\qua\receiveddate\ifx\reviseddate\relax
\else\qquad Revised:\qua\reviseddate\fi\par}} 

\def\lognumber#1{\def\thelognumber{#1}}
\def\volumenumber#1{\def\thevolumenumber{#1}}
\def\volumeyear#1{\def\thevolumeyear{#1}}
\def\papernumber#1{\def\thepapernumber{#1}}
\def\pagenumbers#1#2{\def\startpage{#1}\def\finishpage{#2}}
\def\published#1{\def\publishdate{#1}}

\def\received#1{\def\receiveddate{#1}}
\def\revised#1{\def\reviseddate{#1}}
\def\accepted#1{\def\accepteddate{#1}}

\long\def\asciiabstract#1{\long\def\theasciiabstract{#1}}

\let\\\par\let\thelognumber\relax\let\thevolumenumber\relax
\let\thepapernumber\relax\let\thevolumeyear\relax\let\startpage\relax
\let\finishpage\relax\let\publishdate\relax\let\receiveddate\relax
\let\reviseddate\relax\let\accepteddate\relax\let\theasciititle\relax
\let\theasciiauthors\relax
\let\theasciiabstract\relax

\let\theasciiemail\relax

\ifplaintex
\font\logobig=cmssbx10 scaled 3836
\font\logomed=cmssbx10 scaled 2557
\else
\font\logobig=cmssbx10 scaled 4200
\font\logomed=cmssbx10 scaled 2800
\fi

\long\def\makeagttitle{   
\count0=\startpage
\agt\hfill      
\hbox to 45truept{\vbox to 0pt{\vglue -13truept{\logomed A\kern -.37em{\logobig 
T}\kern -.38em G}\vss}\hss}
\break
{\small Volume \thevolumenumber\ (\thevolumeyear)
\startpage--\finishpage\nl
Published: \publishdate}

\vglue .25truein

{\parskip=0pt\leftskip 0pt plus
1fil\def\\{\par\smallskip}{\Large\bf\thetitle}\par\medskip} \vglue
0.05truein

{\parskip=0pt\leftskip 0pt plus 1fil\def\\{\par}{\sc\theauthors}
\par\medskip}%
 
\vglue 0.03truein 

{\small\leftskip 25truept\rightskip 25truept{\bf Abstract}\stdspace\theabstract

{\bf AMS Classification}\stdspace\theprimaryclass
\ifx\thesecondaryclass\relax\else; \thesecondaryclass\fi\par
{\bf Keywords}\stdspace \thekeywords\par}\vglue 7truept

}   

\ifplaintex
\font\phead=cmsl9 scaled 950
\font\pnum=cmbx10 scaled 913
\font\pfoot=cmsl9 scaled 950
\headline{\vbox to 0pt{\vskip -4.5mm\line{\small\phead\ifnum
\count0=\startpage ISSN 1472-2739 (on-line) 1472-2747 (printed)
\hfill {\pnum\folio}\else\ifodd\count0\def\\{ }%
\ifx\theshorttitle\relax\thetitle\else\theshorttitle\fi\hfill{\pnum\folio}
\else\def\\{ and }{\pnum\folio}\hfill\ifx\theshortauthors\relax\theauthors
\else\theshortauthors\fi\fi\fi}\vss}}
\footline{\vbox to 0pt{\vglue 0mm\line{\small\pfoot\ifnum\count0=\startpage
\copyright\ \gtp\hfill\else
\agt, Volume \thevolumenumber\ (\thevolumeyear)\hfill\fi}\vss}}
\else
\font=cmsl9 scaled 1050
\font\lnum=cmbx10 
\font=cmsl9 scaled 1050
\makeatletter
\def\@oddhead{{\small\lhead\ifnum\count0=\startpage ISSN 1472-2739 
(on-line) 1472-2747 (printed)\hfill {\lnum\number\count0}\else\ifodd\count0
\def\\{ }\ifx\theshorttitle\relax \thetitle \else\theshorttitle\fi\hfill
{\lnum\number\count0}\else\def\\{ and }{\lnum\number\count0}
\hfill\ifx\theshortauthors\relax 
\theauthors\else\theshortauthors\fi\fi\fi}}\def\@evenhead{\@oddhead}
\def\@oddfoot{\small\lfoot\ifnum\count0=\startpage\copyright\ \gtp\hfill\else
\agt, Volume \thevolumenumber\ (\thevolumeyear)\hfill\fi}
\def\@evenfoot{\@oddfoot}
\makeatother
\fi
\let\maketitlepage\makeagttitle
\let\makeshorttitle\maketitlepage
\let\maketitle\maketitlepage

\newwrite\gtoutfile
\long\gdef\makeheadfile{  
{\def\\{, }\def\s{ }
\immediate\openout\gtoutfile head.xxx
\immediate\write\gtoutfile{Proxy-for: \ifx\theasciiauthors\relax
\theauthors\else\theasciiauthors\fi\s<\ifx\theasciiemail\relax\theemail\else\theasciiemail\fi>}
\immediate\write\gtoutfile{\noexpand\\}
\immediate\write\gtoutfile{Authors: \ifx\theasciiauthors\relax
\theauthors\else\theasciiauthors\fi}
{\def\\{ }\immediate\write\gtoutfile{Title: \ifx\theasciititle\relax
\thetitle\else\theasciititle\fi}}
\immediate\write\gtoutfile{Subj-class: GT or SG, GR etc}
\immediate\write\gtoutfile{MSC-class: \theprimaryclass\ifx\thesecondaryclass\relax\else, \thesecondaryclass\fi}
\immediate\write\gtoutfile{Journal-ref: Algebr. Geom. Topol. \thevolumenumber\s
(\thevolumeyear) \startpage-\finishpage}
\immediate\write\gtoutfile{Comments: Published by Algebraic and
Geometric Topology at}
\immediate\write\gtoutfile{\s\s\s  http://www.maths.warwick.ac.uk/agt/AGTVol\thevolumenumber/agt-\thevolumenumber-\thepapernumber.abs.html}
\immediate\write\gtoutfile{\noexpand\\}
\immediate\write\gtoutfile{}
\ifx\theasciiabstract\relax
\immediate\write\gtoutfile{\theabstract}\else
\immediate\write\gtoutfile{\theasciiabstract}\fi
\immediate\write\gtoutfile{}
\immediate\write\gtoutfile{\noexpand\\}
\immediate\write\gtoutfile{}
\immediate\closeout\gtoutfile}}  

\def\maketitlepage{\makeagttitle\makeheadfile}
\let\makeshorttitle\maketitlepage
\let\maketitle\maketitlepage

\lognumber{43}
\volumenumber{4}
\volumeyear{2004}
\papernumber{43}
\pagenumbers{1013}{1040}
\received{5 May 2002} 
\revised{2 July 2003}
\accepted{4 September 2003}
\published{3 November 2004}

\usepackage{ifthen}
\usepackage{amsfonts}
\usepackage{amssymb}
\usepackage{amsopn}
\usepackage{amsmath,amscd}
\usepackage{graphicx}

\theoremstyle{plain}
\newtheorem{thm}{Theorem}[section]
\newtheorem{prop}[thm]{Proposition}

\newtheorem{corol}[thm]{Corollary}

\newtheorem{lem}[thm]{Lemma}

\theoremstyle{definition}
\newtheorem{defn}[thm]{Definition}

\newcommand{\refe}[1]{\eqref{e:#1}}
\newcommand{\reft}[1]{\ref{t:#1}}

\newcommand{\refs}[1]{\ref{s:#1}}

\newcommand{\h}[3]{#1\colon #2\longrightarrow #3}

\begin{document}

\title[The conjugacy problem ]
{The conjugacy problem for relatively\\hyperbolic groups}
\author {Inna Bumagin}

\address{Department of Mathematics and Statistics, Carleton University\\
1125 Colonel By Drive, Herzberg Building\\Ottawa, Ontario, Canada
K1S 5B6} \email{bumagin@math.carleton.ca}

\begin{abstract} Solvability of the conjugacy problem for relatively
hyperbolic groups was announced by Gromov \cite{Gr}. Using the
definition of Farb of a relatively hyperbolic group in the strong
sense~\cite{Farb}, we prove this assertion. We conclude that the
conjugacy problem is solvable for fundamental groups of complete,
finite-volume, negatively curved manifolds, and for finitely
generated fully residually free groups.
\end{abstract}
\asciiabstract{Solvability of the conjugacy problem for relatively
hyperbolic groups was announced by Gromov [Hyperbolic groups, MSRI
publications 8 (1987)]. Using the definition of Farb of a relatively
hyperbolic group in the strong sense [B Farb, Relatively hyperbolic
groups, Geom. Func. Anal. 8 (1998) 810-840], we prove this
assertion. We conclude that the conjugacy problem is solvable for
fundamental groups of complete, finite-volume, negatively curved
manifolds, and for finitely generated fully residually free groups.}

\maketitle

\section{Introduction}
Relatively hyperbolic groups introduced by Gromov \cite{Gr} are
coarsely negatively curved relatively to certain subgroups, called
\emph{parabolic subgroups}. The motivating examples are
fundamental groups of negatively curved manifolds with cusps that
are hyperbolic relative to the fundamental groups of the cusps.
Farb gave his own definition of a relatively hyperbolic group,
using Cayley graphs~\cite[Section 3.1]{Farb} (cf.\ Definition
\reft{rh} below). It was first observed by Szczepanski~\cite{Sz}
that there are groups that satisfy the Farb definition and do not
satisfy the Gromov definition: $\mathbb{Z}\times\mathbb{Z}$ is an
example. For this reason, groups satisfying the Farb definition
are called \emph{weakly relatively hyperbolic}; this terminology
was suggested by Bowditch~\cite{Bow}. Using relative
hyperbolization, Szczepanski~\cite{Sz2} obtained more examples of
weakly relatively hyperbolic groups. Kapovich and Schupp~\cite{KS}
proved that certain Artin groups are weakly relatively hyperbolic.
A weakly relatively hyperbolic group does not have to possess any
nice properties. Osin~\cite{Osin2} showed that there are weakly
relatively hyperbolic groups that are not finitely presentable. He
also constructed an example of a finitely presented weakly
relatively hyperbolic group with unsolvable word problem.

Farb also defined and actually dealt in \cite{Farb} with a
somewhat restricted class of groups, namely, weakly relatively
hyperbolic groups that satisfy the Bounded Coset Penetration (BCP)
property \cite[Section 3.3]{Farb} (cf.\ Definition \reft{BCPp}
below). To prove solvability of the conjugacy problem, we use this
Farb's definition of \emph{relative hyperbolicity in the strong
sense} (see Definition \reft{strong} below).
\begin{thm}\label{t:conj} Let $G$ be a group hyperbolic relative
to a subgroup $H$, in the strong sense. The conjugacy problem is
solvable in $G$, provided that it is solvable in $H$.
\end{thm}
We would like to emphasize importance of the BCP property for
solvability of the conjugacy problem. Collins and Miller \cite{Mi}
give an example of an infinite group $G$ with a subgroup $H$ of
index two (which is therefore, normal in $G),$ so that the conjugacy
problem is solvable in $H$ but is unsolvable in the whole $G.$ In
this example, $G$ is weakly hyperbolic relative to $H$, but
normality of $H$ in $G$ violates the BCP property for the pair
$(G,H)$.

Bowditch~\cite{Bow} elaborated the definitions given by Gromov and
by Farb, and proved that Gromov's definition is equivalent to
Farb's definition of relative hyperbolicity in the strong sense. A
simple alternate proof of the implication ``Gromov's definition
$\Rightarrow$ Farb's definition in the strong sense" can be
derived from the results proved in~\cite{Sz} and~\cite{Bu}
(see~\cite{Bu} for the relevant discussion). It is worth
mentioning that yet another definition of a relatively hyperbolic
group was introduced by Juhasz \cite{J}.

Unlike weakly relatively hyperbolic groups, relatively hyperbolic
groups in the strong sense (which we abbreviate to \emph{relatively
hyperbolic groups}) share many nice properties with word hyperbolic
groups, provided that parabolic subgroups have similar properties.
For instance, Farb proved that the word problem for a relatively
hyperbolic group has ``relatively fast" solution.

\begin{thm}\label{t:F37}{\rm\cite[Theorem 3.7]{Farb}}\qua Suppose $G$ is
strongly hyperbolic relative to a subgroup $H,$  and $H$ has word
problem solvable in time $O(f(n)).$ Then there is an algorithm that
gives an $O(f(n)\log n)$-time solution to the word problem for $G.$
\end{thm}

Arguments that Farb used to prove this latter theorem, imply that
$G$ is finitely presented, if $H$ is; moreover, $G$ has a relative
Dehn presentation. Detailed proofs of these assertions, and of
other basic properties of relatively hyperbolic groups were given
by Osin~\cite{Osin}. Deep results concerning boundaries and
splittings of relatively hyperbolic groups were obtained by
Bowditch~\cite{Bow0},~\cite{Bow1},~\cite{Bow2}. Goldfarb~\cite{Go}
proved Novikov conjectures for these groups, and produced a large
family of relatively hyperbolic groups using strong relative
hyperbolization. Another large family of relatively hyperbolic
groups was produced by Hruska~\cite{Hr}; these are groups acting
properly discontinuously and cocompactly by isometries, on
piecewise Euclidean CAT(0) 2-complexes with isolated flats
property. A topological criterion for a group being relatively
hyperbolic, was obtained by Yaman~\cite{Y}. Dahmani~\cite{Dacl}
proved that a relatively hyperbolic group has a finite classifying
space, if its peripheral subgroups have a finite classifying
space. Rebbechi~\cite{Reb} has shown that relatively hyperbolic
groups are biautomatic, if its peripheral subgroups are
biautomatic. Masur and Minsky \cite{MaM} proved solvability of the
conjugacy problem for mapping class groups, using the fact that a
mapping class group of a surface is weakly hyperbolic relative to
its subgroup that fixes a particular curve on this surface, and
the pair satisfies some additional condition.

In Section~\refs{manifolds} we apply Theorem \reft{conj} to prove
the following.
\begin{thm}\label{t:mnfld} Let $M$ be a complete Riemannian
manifold of finite volume, with pinched negative sectional curvature
and with several cusps. Let $G=\pi_1(M)$ be the fundamental group of
$M.$ Then there is an explicit algorithm to solve the conjugacy
problem for $G.$
\end{thm}
Dehn \cite{Dehn} proved that conjugacy problem for surface groups is
solvable. Cannon \cite{Can} generalized Dehn's proof to all
fundamental groups of closed hyperbolic manifolds. In fact, Cannon's
proof works for the fundamental groups of closed negatively curved
manifolds. Our result can be viewed as a generalization of Cannon's
theorem to the finite volume, noncompact case.

Another example of relatively hyperbolic groups are finitely
generated (f.~g.) fully residually free groups which play an
important role in algebraic geometry over free groups.
\begin{defn} \cite{Ba}\qua A group $L$
is {\it fully residually free}, if for any finite number $n$ of
non-trivial elements $g_1,\dots,$ $g_n$ of $L$ there is a
homomorphism $\varphi$ from $L$ into a free group $F$ so that
$\varphi(g_1),\dots,\varphi(g_n)$ are non-trivial elements of $F.$
\end{defn}
Fully residually free groups are known to have many nice
properties. For this discussion we refer the reader to deep works
of Kharlampovich and Myasnikov \cite{KM}, and also of Sela
\cite{Sd} who introduces the notion of a {\em limit group} and
shows that the classes of limit groups and of f.~g.~fully
residually free groups coincide. The following result is a
conjecture of Sela, proved by Dahmani \cite{Daco}.
Alibegovic~\cite{Al} gave an alternate proof of this conjecture.

\begin{thm}[Dahmani, Alibegovic]\label{t:Daco} Finitely
generated fully residually free groups are relatively hyperbolic
with peripheral structure that consists of the set of their maximal
Abelian non-cyclic subgroups.
\end{thm}

As an immediate corollary of Theorem \reft{Daco} and Theorem
\reft{conj}, we have solvability of the conjugacy problem for f.g.\
fully residually free groups.
\begin{thm}\label{t:cpfrf} The conjugacy problem for finitely
generated fully residually free groups is solvable.
\end{thm}
An alternate proof of this latter assertion, based on using length
functions on fully residually free groups~\cite{MRS}, was given by
Kharlampovich, Myasnikov, Remeslennikov and Serbin in the recent
paper~\cite{KMRS}.

\section{Relatively hyperbolic groups by Farb}\label{s:prelim}

\begin{defn}\cite{Farb}\qua
(Weakly relatively hyperbolic group)\qua\label{t:rh}
Let $G$ be a f.g.\ group, and let $H$ be a f.g.\ subgroup of $G.$ Fix
a set $S$ of generators of $G.$ In the Cayley graph $\Gamma(G,S)$
add a vertex $v(gH)$ for each left coset $gH$ of $H,$ and connect
$v(gH)$ with each $x\in gH$ by an edge of length $\frac 12.$ The
obtained graph $\hat{\Gamma}$ is called a coned-off graph of $G$
with respect to $H.$ The group $G$ is {\it weakly hyperbolic
relative to H} if $\hat{\Gamma}$ is a hyperbolic metric space.
\end{defn}
The above definition depends on the choice of a generating set for
$G.$ Nevertheless, the property of $G$ being weakly hyperbolic
relative to a subgroup $H$ is independent of this choice
\cite[Corollary 3.2 ]{Farb}. Let $u$ be a path in $\Gamma,$ we
define a {\it projection $\hat{u}$ of $u$ into} $\hat{\Gamma}$ in
the special case, when the generating set for $G$ contains a
generating set for $H.$ Reading $u$ from left to right, search for a
maximal subword $z$ of generators of $H.$ If $z$ goes from $g$ to
$g\cdot\bar{z}$ in $\Gamma,$ then we replace the path given by $z$
with the path of length $1$ that goes from the vertex $g$ to the
vertex $g\cdot\bar{z}$ via the cone point $v(gH).$ Do this for each
maximal subword $z$ as above. In general case, projection can be
defined in a similar way: we replace the path given by an element of
$H,$ with a path of length $1$ (see \cite[Section 3.3]{Farb} for
details). We say that $u$ (or $\hat{u}$) travels $\Gamma$-distance
$d_{\Gamma}(g,g\cdot\bar{z})$ in $gH.$ In what follows, we assume
that every path given by a maximal subword $z$ of generators of $H,$
is an $H$-geodesic, in other words we always assume that $z$ is a
path of the shortest $\Gamma$-length that connects $g$ and
$g\cdot\bar{z}.$ Having defined projection $\h
{\rho}{\Gamma}{\hat{\Gamma}}$ by $\rho(u)=\hat{u},$ we can define
relative (quasi)geodesics. Recall that a path $\hat{u}$ with no
self-intersections in $\hat{\Gamma}$ is a $P$-{\em quasi-geodesic}
if for each two points $x,y\in\hat{u}$ the following inequality
holds: $\frac 1P d_{\hat{\Gamma}}(x,y)\leq l_{\hat{u}}(x,y)\leq
Pd_{\hat{\Gamma}}(x,y)$ where $l_{\hat{u}}(x,y)$ denotes the length
of the arc of $\hat{u}$ connecting $x$ and $y$.
\begin{defn} \cite{Farb}\qua({Relative (quasi)geodesics})\qua
\label{t:rqg} If $\hat{u}$ is a geodesic in $\hat{\Gamma},$ then
$u$ is called a {\it relative geodesic} in $\Gamma.$ If $\hat{u}$
is a $P$-quasi-geodesic in $\hat{\Gamma},$ then $u$ is a {\it
relative $P$-quasi-geodesic} in $\Gamma.$
\end{defn}
If $\hat{u}$ passes through some cone point $v(gH),$ we say that
$u$ {\it penetrates} $gH.$ A path $u$ (or $\hat{u}$) is said to be
a {\it path without backtracking} if for every coset $gH$ which
$u$ penetrates, $u$ never returns to $gH$ after leaving $gH.$
\begin{defn}\label{t:BCPp}\cite{Farb}\qua({Bounded Coset Penetration
property})\qua Let a group $G$ be weakly hyperbolic relative to a f.g.\
subgroup $H.$ The pair $(G,H)$ is said to satisfy the {\it Bounded
Coset Penetration (BCP)} property if $\forall P\ge 1,$ there is a
constant $c=c(P)$ so that for every pair $u,v$ of relative
$P$-quasi-geodesics without backtracking, with same endpoints, the
following conditions hold:
\begin{enumerate}
\item \label{t:bcp1} If $u$ penetrates a coset $gH$ and $v$ does
not penetrate $gH,$ then $u$ travels a $\Gamma$-distance of at most
$c$ in $gH.$ \item \label{t:bcp2} If both $u$ and $v$ penetrate a
coset $gH,$ then the vertices in $\Gamma$ at which $u$ and $v$ first
enter (last exit) $gH$ lie a $\Gamma$-distance of at most $c$ from
each other.
\end{enumerate}
\end{defn}
\begin{defn}\label{t:strong}({Strong relative hyperbolicity})\qua
Let $G$ be a f.g.\ group, and let $H$ be a f.g.\ subgroup of $G$. We
say that $G$ is \emph{hyperbolic relative to $H$ in the strong
sense}, if $G$ is weakly hyperbolic relative to $H$, and the pair
$(G,H)$ satisfies the BCP property.
\end{defn}

\section{Notation}\label{s:notation}
Whenever $w$ is a path in $\Gamma$, the projection of $w$ into
$\hat{\Gamma}$ is denoted by $\hat{w}$. Given elements $u$ and $v$
in $G,$ we assume that the equality
\begin{equation} \label{e:uv}
u=gvg^{-1}
\end{equation}
holds for some $g\in G$. We denote by $w$ the closed path in
$\Gamma$ labelled by $ugv^{-1}g^{-1}$. Let $w_u$ and $w_v$ be the
subpaths of $w$ labelled by $u$ and $v^{-1}$, respectively. We
denote by $p$ and $q$ the other two subpaths of $w$. We fix an
orientation of these paths according to the equality
$w=w_upw_vq^{-1},$ so that both paths $p$ and $q$ are labelled by
$g$. Due to the following lemma, we can always assume that $w_u,$
$w_v,$ $p$ and $q$ are relative geodesics.
\begin{lem} \label{t:ugeod} Given an element $x\in G$, one can
find effectively a relative geodesic $\gamma$ that represents $x$.
\end{lem}
\begin{proof} We denote by $x$ both the given element and a $\Gamma$-path that
represents it. Observe that $\hat{x}$ and $\hat{\gamma}$ form a
pair of $l_{\Gamma}(x)$-quasi-geodesics with same endpoints. If
$x$ and $\gamma$ never penetrate the same coset, then
$l_{\Gamma}(\gamma)\leq l_{\Gamma}(x)c(l_{\Gamma}(x))$. If $x$ and
$\gamma$ penetrate a coset $fH$ and $x$ travels along $h_x$ in
$fH$, then $\gamma$ travels a distance bounded by
$l_{\Gamma}(h_x)+2c(l_{\Gamma}(x))$ inside $fH$. Altogether, the
$\Gamma$-length of $\gamma$ is bounded as follows:
$l_{\Gamma}(\gamma)\leq l_{\Gamma}(x)(2c(l_{\Gamma}(x))+1)$. There
are only finitely many elements of $G$ whose length is bounded as
above. Find those that are equal to $x$ in $G$ and take one whose
relative length is minimal possible.
\end{proof}
\begin{corol} \label{t:membH} Given an element $u\in G$, one can
determine effectively whether or not $u$ is in $H$.
\end{corol}
For a fixed set $S$ of generators of $G$, let $S^{\pm 1}$ denote
the set of those generators and their inverses. A product
$g_{i_1}g_{i_2}\dots g_{i_k}$ of elements of $S^{\pm 1}$ is a {\em
reduced word}, if $g_{i_{j+1}}\neq g_{i_j}^{-1}$ for all
$j=1,2,\dots,k-1$. We assume that relative geodesics are labelled
by reduced words. A reduced word $g_{i_1}g_{i_2}\dots g_{i_k}$ is
{\em cyclically reduced}, if $g_{i_1}\neq g_{i_k}^{-1}$. We say
that an element $x\in G$ is {\em cyclically reduced}, if the label
of each relative geodesic $\gamma_x$ that represents $x$ (see
Lemma~\reft{ugeod}), is a cyclically reduced word. Observe that if
$x$ is not cyclically reduced, then a relative geodesic $\gamma_x$
has a proper subpath $\gamma_y$ which is labelled by a cyclically
reduced word. $\gamma_y$ represents an element $y$ of $G$ which is
conjugate to $x$. Since there are only finitely many candidates
for $\gamma_x$ and hence for $\gamma_y$, we can assume that
$\gamma_y$ has the minimal possible relative length. For the
conjugacy problem, we can work with $y$ instead of $x$. If
$\gamma_y$ is a relative geodesic, then we are done. Otherwise, we
will proceed with elements of shorter relative length; therefore,
the process will eventually stop. In what follows, we assume that
$u$ and $v$ are cyclically reduced elements of $G$.

Let $Q=\max\{l_{\Gamma}(u),l_{\Gamma}(v)\}$ denote the maximal
length of $u$ and $v$, and let
$\hat{Q}=\max\{l_{\hat\Gamma}(\hat{u}),l_{\hat\Gamma}(\hat{v})\}$
denote the maximal relative length of $\hat{u}$ and $\hat{v}$. Let
$L=l_{\hat{\Gamma}}(\hat{w})$ be the relative length of $w$, and
let $C=c(L)$ be the constant introduced in Definition \reft{BCPp}.
Observe that the closed path $\hat{w}$ is the concatenation of two
$L$-quasi-geodesic paths as follows: $\lambda_1=w_up$ and
$\lambda_2=w_vq^{-1}$. Denote by $l_H^u$ (or $l_H^v$) the maximal
distance which $u$ (or $v$) travels in an $H$-coset.

\section{Conjugacy problem for hyperbolic groups}\label{s:hyperbolic}
In this section we show that the conjugacy problem for hyperbolic
groups is solvable (see also \cite{Gr},\cite{Lys}). Our proof for
relatively hyperbolic groups is based in part on the extension of
similar ideas to a more general situation, and uses some of the
results proven in this section. To show solvability of the
conjugacy problem for hyperbolic groups, we study properties of
quasi-geodesics in a hyperbolic space. Observe that if $G$ is a
hyperbolic group, then $G$ is hyperbolic relative to the trivial
subgroup so that the coned-off graph $\hat{\Gamma}$ and the Cayley
graph $\Gamma$ of $G$ coincide.
\begin{lem}[Concatenation of two paths
in a geodesic space]\label{t:C} Let $\alpha=\alpha_1\cdot\alpha_2$ be the
concatenation of a geodesic $\alpha_1$ and of a non-empty path
$\alpha_2$ in a geodesic metric space $\Delta$, so that $\alpha$
does not intersect itself and $l_{\Delta}(\alpha_1)\geq
2l_{\Delta}(\alpha_2)$. Then $\alpha$ is a
$(2l_{\Delta}(\alpha_2)+1)$-quasi-geodesic.
\end{lem}
\begin{proof}Let $\beta$ be a geodesic with the same endpoints as
$\alpha$. Then $l_{\Delta}(\alpha_2)\leq\frac 12
l_{\Delta}(\alpha_1)\leq l_{\Delta}(\beta)$. Hence,
$\frac{l_{\Delta}(\alpha_1)+l_{\Delta}(\alpha_2)}{l_{\Delta}(\beta)}
\leq 3$. Now, let $x_i\in\alpha_i$ be a point, let
$\tilde{\alpha}$ be the subpath of $\alpha$ between $x_1$ and
$x_2$, and let $\beta$ be a geodesic joining $x_1$ and $x_2$.
Since $\alpha$ does not intersect itself, $l_{\Delta}(\beta)\geq
1$. It can be readily seen that the maximum possible value of the
ratio $\frac{l_{\Delta}(\tilde{\alpha})}{l_{\Delta}(\beta)}$
equals $2l_{\Delta}(\alpha_2)+1$. Also, note that
$2l_{\Delta}(\alpha_2)+1\geq 3.$
\end{proof}
\begin{lem}\label{t:lemmareL1} If $l_{\hat\Gamma}(\hat{g})\ge
3\hat{Q},$ then the closed path $\hat{w}$ is the concatenation of
two $(2\hat{Q}+1)$-quasi-geodesics. Moreover, without loss of
generality one can assume that the paths
$\lambda_1=\hat{w}_u\hat{p}$ and $\lambda_2=\hat{q}\hat{w}_v^{-1}$
are $(2\hat{Q}+1)$-quasi-geodesics.
\end{lem}
\begin{proof} We prove the assertion for $\lambda_1$, the
proof for $\lambda_2$ is similar. We only need to show that
$\lambda_1$ has no self-intersection.  Assume, $\lambda_1$
intersects itself, which means that $\hat{w}_u$ and $\hat{p}$ have
at least two points in common. Let $x$ be the point were
$\hat{w}_u^{-1}$ and $\hat{p}$ last intersect:
$x=\hat{p}(t)=\hat{w}_u^{-1}(t)$, for some $t$. The path $\hat{w}$
will remain closed if we choose $g$ so that $\hat{p}$ coincides
with $\hat{w}_u^{-1}$ till $x$. Since $u$ is cyclically reduced,
$\hat{w}_u$ and $\hat{q}$ do not intersect. Therefore, the path
$\tilde{\lambda}$ that starts at $\hat{q}(t)$, goes through
$\hat{w}_u$ till $x$ and then through the rest of $\hat{p}$,
satisfies the conditions of Lemma~\reft{C}, so that the first
assertion of the lemma follows. To prove the second assertion,
note that the initial segment $[\hat{q}(t),\hat{w}_u^{-1}(t)]$ of
$\tilde{\lambda}$ represents a cyclic conjugate $\tilde{u}$ of
$u$. Also note that in the closed path $\tilde{w}$ formed by
``cutting off" the common segment of $\hat{w}_u$ and $\hat{p}$,
the path $\tilde{\lambda}$ plays the role of $\lambda_1$.
\end{proof}
\begin{corol}\label{t:corlemmareL1} Let $l_{\hat\Gamma}(\hat{g})\ge
3\hat{Q}$. If $\lambda_1=\hat{w}_u\hat{p}$ backtracks, so that it
can be shortened, then this shorter path $\tilde{\lambda}_1$ is a
$(2\hat{Q}+1)$-quasi-geodesic.
\end{corol}
\begin{proof} $\tilde{\lambda}_1=\alpha_1\cdot\alpha_2$ is
the concatenation of a path $\alpha_1$ with
$l_{\hat\Gamma}(\alpha_1)\le \hat{Q}$ and of a geodesic $\alpha_2$
with $l_{\hat\Gamma}(\alpha_2)\ge 2\hat{Q}$. The assertion follows
from Lemma~\reft{C}.
\end{proof}
Observe that the proofs of Lemma~\reft{C}, Lemma~\reft{lemmareL1}
and Corollary~\reft{corlemmareL1} do not use the assumption that
$\hat{\Gamma}$ is hyperbolic, so that these statements hold for
any geodesic metric space. However, we cannot drop the assumption
that $\hat{\Gamma}$ is hyperbolic, in Corollary~\reft{reL1} and
Lemma~\reft{reL2}  below.
\begin{corol}\label{t:reL1} The paths
$\lambda_1=\hat{w}_u\hat{p}$ and $\lambda_2=\hat{q}\hat{w}_v^{-1}$
stay a bounded distance $K$ from each other in $\hat{\Gamma}$;
moreover, $K$ does not depend on the $\hat{\Gamma}$-length of
$\hat{p}$ (or $\hat{q}$).
\end{corol}
\begin{proof} If $l_{\hat\Gamma}(\hat{g})\le 3\hat{Q},$ then the relative
length of $\lambda_1$ and of $\lambda_2$ is bounded by $4\hat{Q}$,
so that $\lambda_1$ and $\lambda_2$ stay a distance bounded by
$2\hat{Q}$ from each other. Otherwise, by Lemma~\reft{lemmareL1},
$\lambda_1$ and $\lambda_2$ are $(2\hat{Q}+1)$-quasi-geodesics
with common endpoints. Therefore, they stay a bounded distance
$N(2\hat{Q}+1)$ from each other. In order to obtain the claim, set
\begin{equation} \label{e:K}
K=\max\{2\hat{Q},N(2\hat{Q}+1)\}.
\end{equation}
Note that the initial point of both $\hat{w}_u$ and $\hat{q}$ is
the identity $1_{\hat{\Gamma}},$ and that the initial point of
$\hat{p}$ coincides with the terminal point of $\hat{w}_u.$
\end{proof}
\begin{lem}\label{t:reL2} Start at the initial points of $\hat{p}$ and
$\hat{q}$ and move along these paths with the unit speed. If
$\hat{p}$ and $\hat{q}$ are long enough paths, then there exist
two numbers $t_1,t_2$ satisfying $1\le t_1<t_2\le
l_{\hat\Gamma}(\hat{p}),$ so that the following condition holds.
For each integer $t$ where $t_1\le t\le t_2,$ the
$\hat{\Gamma}$-length of a geodesic path $\hat\gamma$ joining
$\hat{p}(t)$ and $\hat{q}(t)$ can be bounded in terms of the
$\hat\Gamma$-lengths of $u$ and $v.$
\end{lem}
\begin{proof} Assume that $\hat{p}$ and $\hat{q}$ are longer than
$3\hat{Q}+2\delta$ so that we can consider values of $t_1,t_2$ as
follows:
\begin{equation} \label{e:t1t2}
l_{\hat\Gamma}(\hat{u})+\delta=t_1<t_2=
l_{\hat\Gamma}(\hat{p})-(l_{\hat\Gamma}(\hat{v})+\delta).
\end{equation}
Denote by $\hat{p'}$ (or $\hat{q'})$ the subpath  of $\hat{p}$ (or
$\hat{q})$ between $\hat{p}(t_1)$ (or $\hat{q}(t_1))$ and
$\hat{p}(t_2)$ (or $\hat{q}(t_2)).$ Because of our choice of $t_1$
and $t_2,$ the geodesics $\hat{p'}$ and $\hat{q'}$ stay the
distance $K$ from each other ($K$ as in \refe{K}). Furthermore,
the distance between the initial points of the paths $\hat{p'}$
and $\hat{q'}$ is bounded above as follows:
$d_{\hat\Gamma}(\hat{p}(t_1),\hat{q}(t_1))\le l_0,$ where
$l_0=3\hat{Q}+2\delta.$

Note that $d_{\hat\Gamma}(\hat{p}(t_1),\hat{p}(t))=
d_{\hat\Gamma}(\hat{q}(t_1),\hat{q}(t))=|t-t_1|.$ Consider
$\hat{\gamma}(t)$ for $t_1<t<t_2.$ There is a point
$x=\hat{q}(t_x)$ so that $d_{\hat\Gamma}(\hat{p}(t),x)\le K.$
Assume, $t_1\le t_x\le t,$ so that
$d_{\hat\Gamma}(\hat{q}(t),x)=|t-t_1|-d_{\hat\Gamma}(\hat{q}(t_1),x).$
Since $\hat{p'}$ is a geodesic, we have that $|t-t_1|\le
l_0+d_{\hat\Gamma}(\hat{q}(t_1),x)+K.$ Hence,
$d_{\hat\Gamma}(\hat{q}(t),x)\le l_0+K.$ If $t_x\ge t,$ then we
have that
$d_{\hat\Gamma}(\hat{q}(t_1),x)=|t-t_1|+d_{\hat\Gamma}(\hat{q}(t),x)\le
K+|t-t_1|+l_0$, so that $d_{\hat\Gamma}(\hat{q}(t),x)\le l_0+K$ as
well. In both cases, we conclude that
\begin{equation} \label{e:gamma}
l_{\hat\Gamma}(\hat{\gamma}(t))=d_{\hat\Gamma}(\hat{p}(t),\hat{q}(t))\le
3\hat{Q}+2\delta+2K,
\end{equation}
which implies the claim.
\end{proof}
As a consequence, we get the following theorem which is the main
result of this section.
\begin{thm} \label{t:hyperbolic} If $G$ is a hyperbolic group,
then the conjugacy problem in $G$ is solvable.
\end{thm}
\begin{proof} It follows immediately from Lemma
\reft{reL2} that for $\hat{p}$ and $\hat{q}$ long enough, one can
find two integers $s_1$ and $s_2$ satisfying the double inequality
$t_1\le s_1<s_2\le t_2$, so that two geodesics $\hat{\gamma}_i$
connecting $\hat{p}(s_i)$ with $\hat{q}(s_i)$ $(i=1,2)$ have the
same $\hat{\Gamma}$-length. Therefore, the compactness of balls of
a given radius in $\hat\Gamma$ implies that for $\hat{p}$ and
$\hat{q}$ long enough, we can find $\hat{\gamma}_1,\hat{\gamma}_2$
as above so that these geodesics represent the same element $x$ of
$G.$ In this case, both $u$ and $v$ are conjugate to $x$ in $G$.
Moreover, we can cut off the segments
$[\hat{p}(s_1),\hat{p}(s_2)]$ and $[\hat{q}(s_1),\hat{q}(s_2)]$ of
$\hat{p}$ and $\hat{q}$, respectively, and obtain a shorter
element conjugating $u$ and $v.$ Thus, if two elements $u$ and $v$
of a hyperbolic group are conjugate, then the minimal possible
length of a conjugating element $g$ is bounded in terms of the
lengths of $u$ and $v$ (cf.\ \cite[Lemma 10]{Lys}). Since bounded
balls in $\hat{\Gamma}$ are compact and word problem in $G$ is
solvable, the assertion follows.
\end{proof}

\section{Relatively hyperbolic groups}\label{s:relhyp}
We assume that the conjugacy problem is solvable in $H.$
Therefore, given $u,v\in H,$ we can determine effectively whether
or not $u$ and $v$ are conjugate in $H$. If this is the case, then
$\hat{w}$ is the null-path in $\hat{\Gamma}$. In what follows, we
assume that $L=l_{\hat\Gamma}(\hat{w})>1$.
\begin{defn}\label{t:closednobacktrack}(Closed path without
backtracking)\qua We will say that a {\it closed path $\hat{w}$ does
not backtrack to a coset} $gH$ which it penetrates, if $\hat{w}$
is the concatenation of two paths $\hat{u}$ and $\hat{v}^{-1}$ so
that $\hat{u}$ is a path without backtracking which penetrates
$gH,$ and $\hat{v}$ does not penetrate $gH.$ We will say that the
{\it closed} path $\hat{w}$ is a {\it path without backtracking},
if for any coset $gH$ which it penetrates, $\hat{w}$ does not
backtrack to $gH.$
\end{defn}
The following results which will be used later on, are
straightforward consequences of the Definition~\reft{BCPp} of the
BCP property. Recall that $C=c(L)$ is the constant introduced in
Definition \reft{BCPp}.
\begin{lem}\label{t:cboup} If $\hat{w}$ penetrates a coset
$gH$ and does not backtrack to it, then $\hat{w}$ travels in $gH$
a $\Gamma$-distance of $C=c(L)$ at most.
\end{lem}
\begin{proof} Since each relative $P$-quasigeodesic is in particular
an $R$-quasigeodesic for $P<R,$ we have that $c(P)\leq c(R).$ Let
$\hat{w}$ travel in $gH$ along an $H$-geodesic $h.$ Consider the
subpath $\hat{w}_1=g_1hg_2$ of $\hat{w}$. The relative length of
$\hat{w}_1$ equals $3$ so that $\hat{w}_1$ is a $3$-quasigeodesic.
If $\hat{w}_1$ is a closed path, then we have found a pair of two
relative $2$-quasigeodesics, one of which penetrates the coset
$gH,$ and the other does not penetrate $gH.$ Therefore, in this
case $L_{\Gamma}(h)\le c(2)\le c(L).$ Now, assume $\hat{w}_1$ is
not a closed path. Let $\hat{w}_2$ be so that
$\hat{w}=\hat{w}_1\hat{w}_2.$ Since $\hat{w}$ does not backtrack
to $gH,$ $\hat{w}_2$ does not penetrate this coset. If $\hat{w}_2$
backtracks to a coset different from $gH,$ then we can shorten
$\hat{w}_2$ each time when backtracking occurs. Indeed, assume
$\hat{w}_2$ leaves a coset $fH$ at some point $x$ and enters this
coset at $y,$ later on. We replace the subpath of $\hat{w}_2$
joining $x$ and $y,$ with an $H$-geodesic $h_{x,y}$ joining these
points. Finally, we obtain a path $\hat{w}'_2$ without
backtracking, with $L_{\hat{\Gamma}}(\hat{w}'_2)\le L-3.$ Hence,
$\hat{w}'_2$ is a relative $(L-3)$-quasigeodesic, which does not
penetrate $gH.$ Therefore, the pair $\hat{w}_1$ and $\hat{w}'_2$
of relative $(L-3)$-quasigeodesics satisfies the first part of
Definition \reft{BCPp}, and it follows that $L_{\Gamma}(h)\le
c(L-3)\le c(L).$
\end{proof}
\begin{figure}[ht!]
\centering
\includegraphics[width=.3\textwidth]{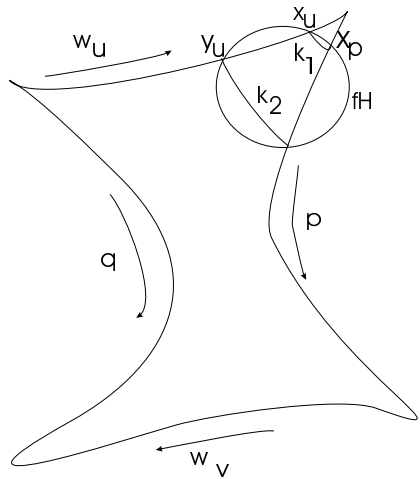}
\caption{Illustrated above is the case when $\hat{w}_u$ and
$\hat{p}$ penetrate a coset $fH$, which $\hat{q}$ and $\hat{w}_v$
do not penetrate (Lemma~\reft{A}). \label{fig:fig53}}
\end{figure}
\begin{lem}\label{t:A} Assume that $\hat{w}_u$ and $\hat{p}$ penetrate a coset
$fH$, but neither $\hat{q}$ nor $\hat{w}_v$ penetrates $fH$. Let
$k_1$ (or $k_2$) be a $\Gamma$-geodesic joining the points where
$\hat{w}_u^{-1}$ and $\hat{p}$ first enter (or last exit) the
coset $fH$. Then $l_{\Gamma}(k_1)\leq c(2)$ and
$l_{\Gamma}(k_2)\leq C$ (Figure~\ref{fig:fig53}).
\end{lem}
\begin{proof} Denote by $x_p$ and $x_u$ (or $y_p$ and $y_u$) the
endpoints of $k_1$ (or of $k_2$). To prove the inequality for
$l_{\Gamma}(k_1)$, consider the following pair of
$2$-quasi-geodesics: $\alpha$ is the concatenation of the initial
subgeodesic of $\hat{p}$ (ending at $x_p$) and $k_1$, $\beta$ is
the initial subgeodesic of $\hat{w}_u^{-1}$ (ending at $x_u$).
Note that $\alpha$ and $\beta$ have common endpoints and do not
backtrack. Therefore, the first part of the definition of the BCP
property implies the claim for $k_1$. To prove the inequality for
$l_{\Gamma}(k_2)$, consider the closed path $\hat{w}'$ which is
the concatenation of the following geodesics: $k_2$, the terminal
subgeodesic of $\hat{w}_u^{-1}$ (starting at $y_u$), $\hat{q}$,
$\hat{w}_v^{-1}$ and the initial subgeodesic of $\hat{p}^{-1}$
(ending at $y_q$). As $\hat{w}'$ is shorter than $\hat{w}$,
Lemma~\reft{cboup} implies the claim for $k_2$.
\end{proof}
\begin{corol} \label{t:corolA} Let $\lambda_1=w_up$ and
$\lambda_2=w_vq^{-1}$ be relative $P$-quasi-geodesics for some
$P>0$. Let $p$ travel along the path $h_g$ in a coset $fH$. If $p$
and $q$ do not penetrate the same coset, then $l_{\Gamma}(h_g)\leq
2Q+c(P)+2c(2).$
\end{corol}
\begin{proof} Let $\hat{w}_u$ and $\hat{p}$ penetrate a coset
$fH$. The proof of Lemma~\reft{A} works verbatim in the case when
$\hat{w}_v$ does not penetrate the coset $fH$, but instead of
Lemma~\reft{cboup} use Corollary~\reft{corlemmareL1}. Now, assume
$\hat{w}_v$ penetrates $fH$ also. Let $k_1$ be as in the statement
of Lemma~\reft{A}, let $k_2$ be a geodesic joining the points at
which $\hat{w}_v$ and $\hat{p}^{-1}$ first enter $fH,$ and let
$k_3$ join the points at which $\hat{w}_u$ and $\hat{w}_v^{-1}$
first enter $fH$. It follows from the argument used to prove
Lemma~\reft{A} that $l_{\Gamma}(k_i)\leq c(2)$ for $i=1,2$.
Moreover, it can be readily seen that the concatenation of the
initial segment of $\hat{w}_u$ and of $k_3$ is a relative
$3$-quasi-geodesic. By assumption, the concatenation of
$\hat{q}^{-1}$ and of the initial segment of $w_v^{-1}$ is a
relative $P$-quasi-geodesic which does not penetrate $fH$, so that
$l_{\Gamma}(k_3)\leq c(P)$, which implies the claim.
\end{proof}
\textbf{Remark.} Corollary~\reft{corolA} holds for each closed
path which is the concatenation of four geodesics. Indeed, the
fact that $\hat{p}$ and $\hat{q}$ have the same label is never
used in the proof of this statement.
\begin{figure}[ht!]\centering
\includegraphics[width=.75\textwidth]{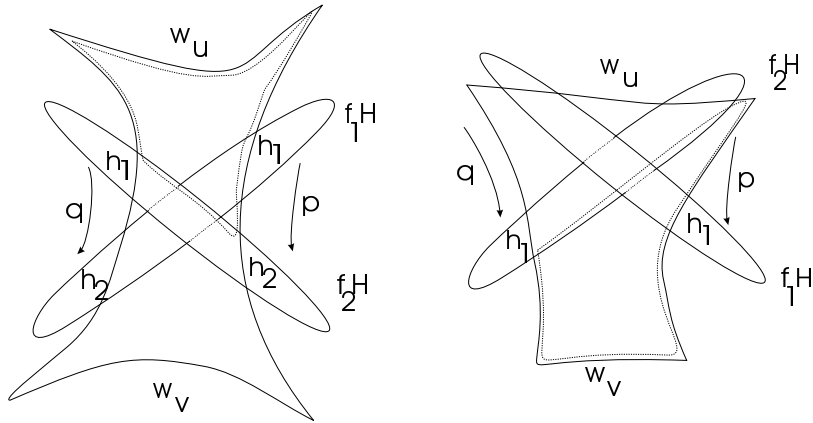}
\caption{Examples of ``skew" cosets (Lemma~\reft{B}). The dotted
line shows a closed path which does not backtrack and travels a
bounded distance in $f_1H$. \label{fig:fig55}}
\end{figure}
\begin{lem}\label{t:B} Assume that $\hat{p}$ and $\hat{q}$ penetrate two ``skew"
cosets $f_1H$ and $f_2H$ so that $\hat{p}$ travels in $f_1H$ along
$h_1$ and $\hat{q}$ travels in $f_2H$ along $h_1$, and either
\begin{enumerate}
 \item \label{e:pqonly} $\hat{w}_u$ and $\hat{w}_v$ penetrate neither of these
 cosets, and $\hat{p}$ travels in $f_2H$ along $h_2$ while $\hat{q}$ travels in $f_1H$ along
 $h_2$ (Figure~\ref{fig:fig55}, left), or
 \item \label{e:pw_uq} $\hat{w}_u$ penetrates first $f_1H$ and then $f_2H$,
 $\hat{w}_v$ penetrates neither of these cosets, $\hat{p}$ does not penetrate
 $f_2H$, and $\hat{q}$ does not penetrate $f_1H$ (Figure~\ref{fig:fig55}, right).
\end{enumerate} Then $l_{\Gamma}(h_1)\leq C$.
\end{lem}
\begin{proof} We prove the assertion only in the case~\refe{pqonly}, the other
case is similar. Consider the closed path $\hat{w}'=\hat{p}_1\circ
k\circ\hat{q}_1^{-1}\circ \hat{w}_u$, where $k$ is an $H$-
geodesic in $f_2H$, and $\hat{p}_1$ (or $\hat{q}_1$) is the
initial segment of $\hat{p}$ (or $\hat{q}$) that terminates at the
point where $\hat{p}$ (or $\hat{q}$) first enters the coset
$f_2H$. This path satisfies the conditions of Lemma~\reft{cboup}
and is shorter than $\hat{w}$; the assertion~\refe{pqonly}
follows.
\end{proof}
\subsection{Cascades}\label{s:cascades}
\begin{figure}[ht!]\centering
\includegraphics[width=.8\textwidth]{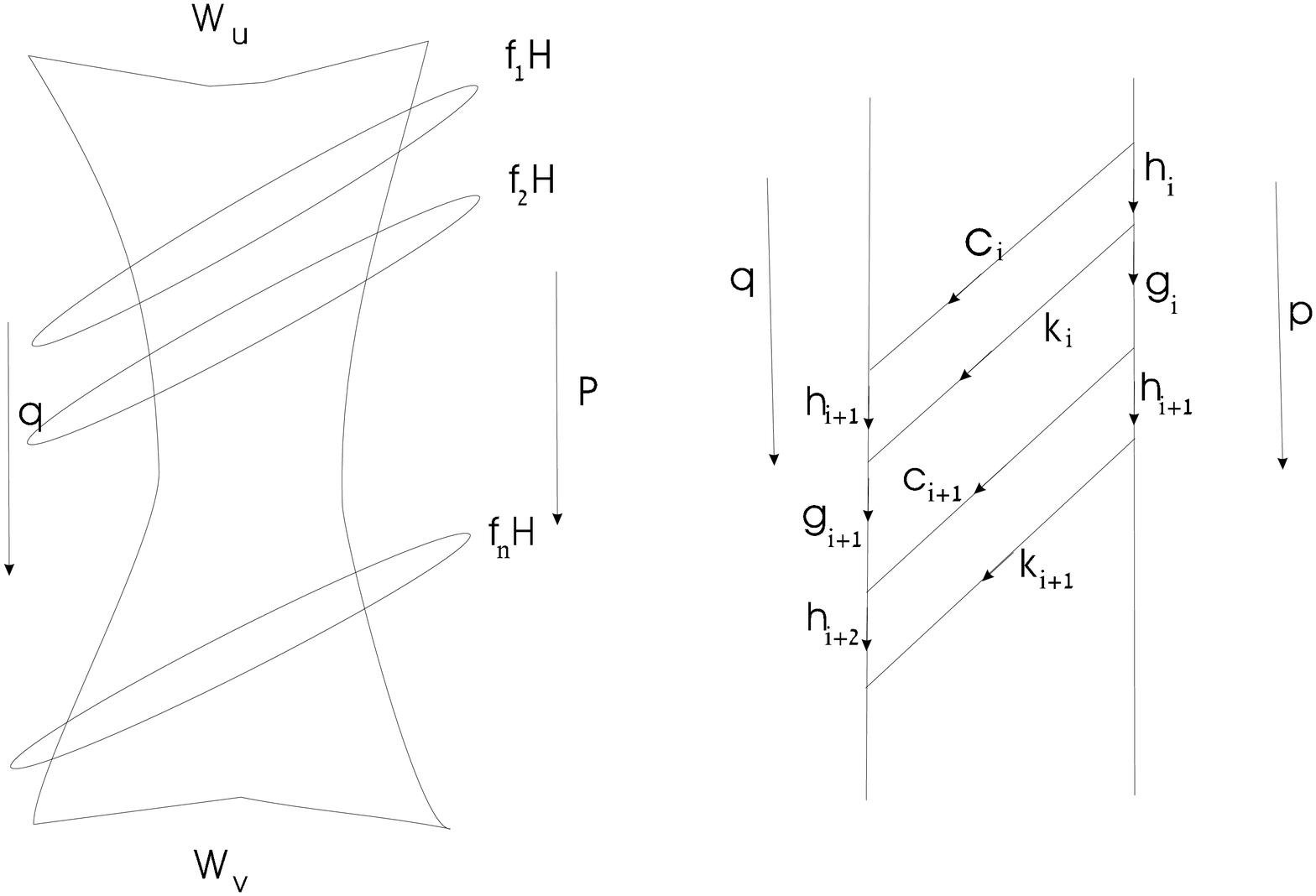}
\caption{Left: Cascade effect. Right: Three consecutive floors of
a cascade, two $H$-floors and a $G$-floor between them.
\label{fig:cascade}}
\end{figure}

So far, we were able to apply directly the definition of the BCP
property in order to bound the $\Gamma$-length of a subpath of $w$
in terms of the $\hat{\Gamma}$-length of $\hat{w}$. It is possible
unless $\hat{p}$ and $\hat{q}$ penetrate the same cosets. Contract
each vertex of $\hat{\Gamma}$ to a point; if $\hat{p}$ and
$\hat{q}$ penetrate the same cosets, then $\hat{w}$ will turn into
a sequence of digons with two triangles at both ends of it. In the
case when these triangles are isosceles, the length of the
conjugating element $g$ cannot be bounded in general. We consider
this case in Section~\refs{conjh} below. If the triangles are not
isosceles, then we have {\em cascade effect} defined as follows.
\begin{defn}\label{t:defcascade}(Cascade effect)\qua We say
that in the path $w=w_upw_vq^{-1}$ {\it cascade effect} occurs if
the following condition holds. There are subwords
$h_1,h_2,\dots,h_{n+1}\in H$ of $g$ and cosets
$f_1H,f_2H,\dots,f_{n}H$, which both $p$ and $q$ penetrate in the
same order, and so that either $p$ travels in the coset $f_iH$
along $h_i$ and $q$ travels in the coset $f_iH$ along $h_{i+1}$
(Figure~\ref{fig:cascade}, left), or $p$ travels in the coset
$f_iH$ along $h_{i+1}$ and $q$ travels in the coset $f_iH$ along
$h_{i},$ for each $i=1,2,\dots,n$. The number $n$ is called the
{\it length of the cascade.}
\end{defn}
The ``tower" corresponding to a cascade is an alternating sequence
of $H$-floors and $G$-floors (Figure~\ref{fig:cascade}, right);
each $H$-{\em floor} is a closed path $h_ik_ih_{i+1}^{-1}c_i^{-1}$
in $H$, each $G$-{\em floor} is a closed path
$g_ic_{i+1}g_{i+1}^{-1}k_i^{-1}$ where $g_i\in G\setminus H$ and
$k_i,c_i\in H$. Each $G$-floor (or $H$-floor) has two neighboring
$H$-floors (or $G$-floors) glued to it along $\Gamma$-geodesic
arcs corresponding to $k_i$ and $c_{i+1}$ (or $c_i$ and $k_i$).
Observe that in a cascade of length $n$, the number of $H$-floors
is $n$ and the number of $G$-floors is $n-1$. It can be readily
seen that the length of the subwords $h_i$ of $g$ can be bounded
in terms of the length of the cascade as follows:
$l_{\Gamma}(h_i)\le l_{\Gamma}(h_1)+2ic(2)\le
l_{\Gamma}(h_1)+2nc(2)$ (see the part~\refe{ki} of
Lemma~\reft{firstlast} below); the main difficulty is to handle
cases when $n$ is large. In Lemma~\reft{cascade} below we obtain a
bound on the length of $h_i$ which does not depend on the length
of a cascade. Lemma~\reft{firstlast} asserts that the
$\Gamma$-length of each $k_i$ and of each $c_i$ can be bounded in
terms of the relative length of $u$ and $v$. Since we can skip
several $H$-floors in the beginning and in the end of a cascade,
we will always assume that neither $w_u$ nor $w_v$ penetrate the
cosets where $H$-floors of a cascade are located. We set
$C_0=c(7Q)$.
\begin{lem}\label{t:firstlast} With the notation above, the length
of the $\Gamma$-geodesics $c_i$ and $k_i$ can be bounded as
follows.
\begin{enumerate}
  \item \label{e:ki} $l_{\Gamma}(k_i)\le c(2)$ for $1\le i<n$, and
$l_{\Gamma}(c_i)\le c(2)$ for $1<i\le n$.
  \item \label{e:k1n} $l_{\Gamma}(k_n)\le C_0$ and
  $l_{\Gamma}(c_1)\le C_0$.
\end{enumerate}
\end{lem}
\begin{proof} As each $G$-floor is the concatenation of two
$2$-geodesics, the assertion~\refe{ki} follows immediately. We
prove the assertion~\refe{k1n} for $c_1$, the proof for $k_n$ is
similar. Consider the closed path $w'$ which is the concatenation
of $w_u$, $p_1$, $c_1$ and $q_1^{-1}$, where $p_1$ and $q_1$ are
the initial segments of $p$ and $q$, respectively.
Without loss of generality, assume that $l_{\Gamma}(p_1)\le
l_{\Gamma}(q_1)$. The concatenation $\alpha=p_1\cdot c_1$ is a
relative $2$-quasi-geodesic.

If $p_1$ and $q_1$ penetrate a coset $fH$, then a subpath of
$\alpha$ as well as the concatenation of a $\Gamma$-geodesic which
travels in $fH$ and of $q'_1$ are relative $2$-quasi-geodesics
with common endpoints. Hence, $l_{\Gamma}(c_1)\le c(2)$ in this
case. In what follows, we assume that $p_1$ and $q_1$ do not
penetrate the same coset.

We distinguish the two cases as follows:
\begin{enumerate}
  \item \label{e:long} $\hat{q}_1$ is ``long": $l_{\hat{\Gamma}}(\hat{q}_1)\ge
  3Q>2l_{\hat{\Gamma}}(\hat{w}_u)$.
  \item \label{e:short} $\hat{q}_1$ is ``short": $l_{\hat{\Gamma}}(\hat{q}_1)<3Q$.
\end{enumerate}
In the case~\refe{long} our argument below shows that
$l_{\Gamma}(c_1)\le c(2Q+1)$. Indeed, the concatenation
$\hat{\beta}=\hat{w}_u^{-1}\cdot \hat{q}_1$ is a
$(2Q+1)$-quasi-geodesic, according to Lemma~\reft{C}. If
$\hat{\beta}$ backtracks, then by Corollary~\reft{corlemmareL1}, a
shorter path $\tilde{\beta}$ without backtracking and same
endpoints as $\hat{\beta}$, is a $(2Q+1)$-quasi-geodesic.
Furthermore, if $\hat{p}_1$ and $\hat{w}_u^{-1}$ penetrate a coset
$fH$, then set $\tilde{\alpha}$ to be the subpath of
$\hat{\alpha}$ which begins at the point where $\hat{p}_1$ exits
$fH$, and adjust $\tilde{\beta}$ accordingly. In any case, we have
a pair of $(2Q+1)$-quasi-geodesics that satisfies the first part
of the definition of the BCP property, which implies the claim in
the case when $\hat{q}_1$ is ``long".

In the case~\refe{short} $l_{\hat{\Gamma}}(\hat{w}')\le 7Q$, and
so by Lemma~\reft{cboup}, $l_{\Gamma}(c_1)\le C_0=c(7Q)$, as
claimed.
\end{proof}
\begin{lem}\label{t:cascade} Assume that in the word
$w$ the cascade effect occurs. Let $p$ travel in the coset $f_iH$
along $h_i$. Then $l_{\Gamma}(h_i)\le l_{\Gamma}(h_1)+2C_0+2c(2)$,
for $i=2,3,\dots, n-1$.
\end{lem}
\begin{proof} Since each $G$-floor is the concatenation of two
$2$-geodesics, $g_i$ and $g_{i+l}$ travel a $\Gamma$-distance
bounded by $c(2)$ in each coset they penetrate. By
Lemma~\reft{firstlast} case~\refe{ki}, $l_{\Gamma}(h_2)\le
l_{\Gamma}(h_1)+2c(2)$ and $l_{\Gamma}(h_3)\le
l_{\Gamma}(h_1)+4c(2)$, so that in what follows, we assume that
$i\geq 4$. To show that the $\Gamma$-length of $h_i$ is bounded,
we glue $i-2$ consecutive $G$-floors of the cascade along $g_j$
where $j=2,3,\dots,i-2$, and $i-1$ consecutive $H$-floors of the
cascade along $h_j$, for $j=2,3,\dots,i-1$. We have the following
equalities:
\begin{align*}
k_1k_2\dots k_{i-2} & =
g_1c_2c_3\dots c_{i-1}g_{i-1}^{-1} \\
k_1k_2\dots k_{i-2} & = h_1^{-1}c_1c_2\dots
c_{i-2}c_{i-1}h_{i}k_{i-1}^{-1}
\end{align*}
Denote $\bar{c}=c_2c_3\dots c_{i-1}$ and
$\bar{h}=h_1^{-1}c_1\bar{c}h_{i}k_{i-1}^{-1}$; it follows from the
equalities above that $g_1\bar{c}g_{i-1}^{-1}=\bar{h}$. Therefore,
$g_1\bar{c}$ and $\bar{h}g_{i-1}$ form a pair of
$2$-quasi-geodesics with common endpoints. Hence, the BCP property
implies that the $\Gamma$-length of both $\bar{c}$ and $\bar{h}$
is bounded by $c(2)$. Since
$h_i=(c_1\bar{c})^{-1}h_1\bar{h}k_{i-1}$, according to
Lemma~\reft{firstlast}, we have that $l_{\Gamma}(h_i)\le
l_{\Gamma}(h_1)+C_0+3c(2)$.
\end{proof}
\subsection{Relatively short conjugating elements are
short}\label{s:short}
\begin{figure}[ht!]\centering
\includegraphics[width=.7\textwidth]{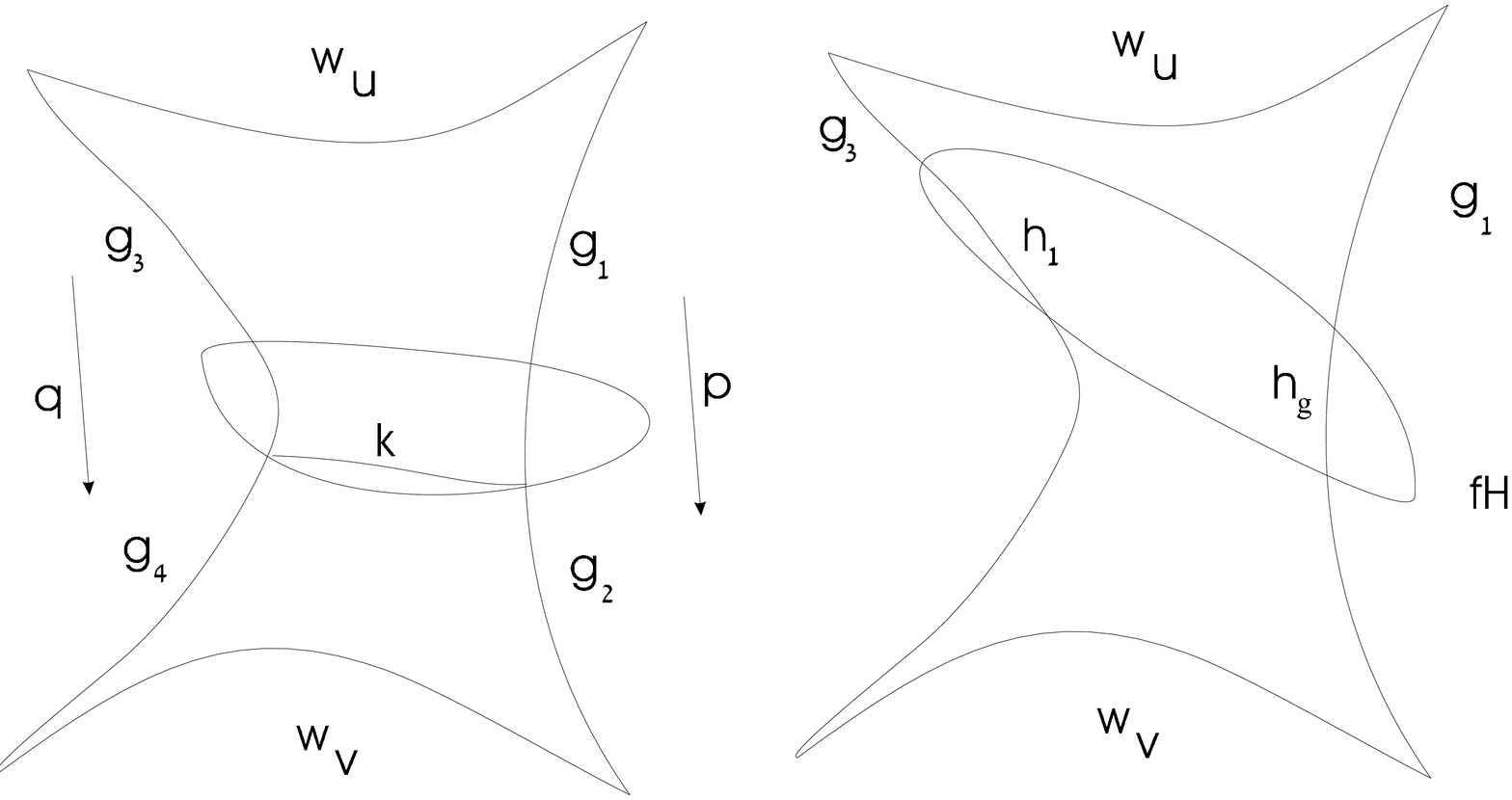}
\caption{Illustrated above is Lemma~\reft{mtl}, case~\refe{two},
when $\hat{p}$ and $\hat{q}$ penetrate a coset $fH.$
\label{fig:fig585}}
\end{figure}
 The following lemma enables one to bound the
$\Gamma$-length of $g$ in terms of the relative length of
$\hat{g}$, $\hat{u}$ and $\hat{v}$.
\begin{lem} \label{t:mtl} Let $\hat{p}$ penetrate a left
coset $fH$ when moving along an $H$-geodesic $h_g$. Assume that
$\hat{w}$ backtracks to $fH.$
 \begin{enumerate}
  \item \label{e:four} If $\hat{p},$ $\hat{w_u},$ $\hat{q}$
  and $\hat{w_v}$ all penetrate the coset $fH,$ then either
    \begin{enumerate}
    \item   \label{e:foura} 
     $u$ and $v$ are conjugate in $G$ to $k\in H$
     with the length bounded as follows:
    $l_{\Gamma}(k)\leq Q+2C,$ or
    \item  \label{e:fourb} 
    $l_{\Gamma}(h_g)\leq Q+2C.$
    \end{enumerate}
  \item \label{e:gug}If $\hat{p},$ $\hat{w_u}$ and $\hat{q}$
  penetrate $fH,$ and $\hat{w_v}$ does not penetrate it, then either
   \begin{enumerate}
    \item  \label{e:guga} 
    $u$ and $v$ are conjugate in $G$ to $k\in H$ with
    $l_{\Gamma}(k)\leq C,$ or
     \item  \label{e:gugb} 
     $l_{\Gamma}(h_g)\leq Q+3C+3c(2).$
    \end{enumerate}
  \item \label{e:guv}If $\hat{p},$ $\hat{w_u}$ and $\hat{w_v}$
  penetrate $fH,$ and $\hat{q}$ does not penetrate it, then
  $l_{\Gamma}(h_g)\leq 2Q+3C.$
  \item \label{e:gu} If $\hat{p}$ and $\hat{w_u}$ penetrate a coset $fH,$ and
  neither $\hat{w_v}$ nor $\hat{q}$ penetrates $fH$, then
  $l_{\Gamma}(h_g)\leq Q+2C.$
  \item \label{e:two} If $\hat{p}$ and $\hat{q}$ penetrate
  $fH,$ and neither $\hat{w_u}$ nor $\hat{w_v}$ penetrates it, then either
     \begin{enumerate}
    \item  \label{e:twoa} 
    $u$ and $v$ are conjugate in $G$ to $k\in H$ with
    $l_{\Gamma}(k)\leq C,$ or
     \item  \label{e:twob} 
     $l_{\Gamma}(h_g)\leq Q+3C+3c(2).$
    \end{enumerate}
 \end{enumerate}
\end{lem}
\begin{figure}[ht!]\centering
\includegraphics[width=.4\textwidth]{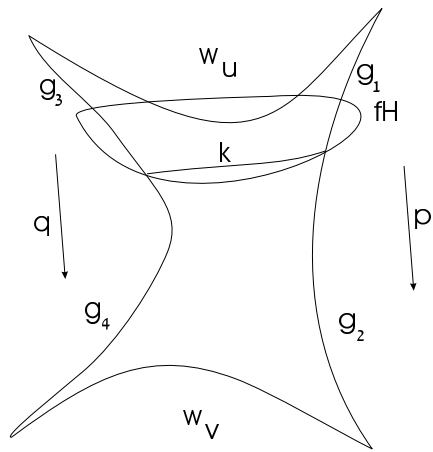}
\caption{Illustrated above is Lemma~\reft{mtl}, case~\refe{gug},
when $\hat{p},$ $\hat{w_u}$ and $\hat{q}$ penetrate a coset $fH,$
and $u$ and $v$ are conjugate to $k$. \label{fig:fig582a}}
\end{figure}
\begin{proof} Case \refe{guv} is a particular case of
Corollary~\reft{corolA}. The statement in the case \refe{gu}
follows from Lemma~\reft{A}. In the proof below, we use the
following notation. Let $g=g_1h_gg_2$ (we denote by $g_i$ the
sub-paths of $p$ and also their labels), and let $fH$ be the coset
in which $\hat{p}$ travels along $h_g$. In the case~\refe{two},
$\hat{q}$ penetrates the coset $fH$ also; denote by $h_f$ the
$H$-geodesic in $fH,$ along which $\hat{q}$ moves there.
Therefore, $g=g_3h_fg_4,$ where $g_3,g_4$ denote both sub-paths of
$q$ and their labels. If
$l_{\hat{\Gamma}}(\hat{g_1})=l_{\hat{\Gamma}}(\hat{g_3})$
(Figure~\ref{fig:fig585}, left), then necessarily $g_1=g_3,$
$h_g=h_f,$ and $g_2=g_4.$ Let $k$ be a $H$-geodesic joining the
points where $p$ and $q$ last exit $fH.$ We have that
$k=h_g^{-1}g_1^{-1}ug_1h_g$, and $k=g_2vg_2^{-1}$. Moreover, the
closed path $g_2w_vg_2^{-1}k$ has length less than $L$ and
satisfies the conditions of Lemma~\reft{cboup}. Hence, we obtain
\refe{twoa}. Now, assume that
$l_{\hat{\Gamma}}(\hat{g_1})>l_{\hat{\Gamma}}(\hat{g_3})$
(Figure~\ref{fig:fig585}, right). Let $k_1$ (or $k_2$) be a
geodesic joining the points where $p$ and $q$ first enter (or last
exit) $fH$; both closed paths (one goes through $k_1$ and $w_u$
and the other one goes through $k_2$ and $w_v$) that we obtain,
are shorter than $w$ and satisfy the conditions of
Lemma~\reft{cboup}. Hence, $|l_{\Gamma}(h_f)-l_{\Gamma}(h_g)|\le
2C.$ Furthermore, $h_f$ is a subword of $g_1.$ Let $f_1H$ be the
coset that $\hat{p}$ penetrates along $h_f.$ If $\hat{w}$ does not
backtrack to $f_1H$, then one can apply Lemma~\reft{cboup} and
obtain that $l_{\Gamma}(h_f)\leq C$. If $\hat{w}$ backtracks to
$f_1H$ but in $\hat{w}$ cascade effect does not occur, then we
apply either Lemma~\reft{B} or Lemma~\reft{A}, and conclude that
$l_{\Gamma}(h_f)\leq l_{\Gamma}(v)+C+c(2)$, so that
$l_{\Gamma}(h_g)\leq l_{\Gamma}(v)+3C+c(2)$. If in $\hat{w}$
cascade effect occurs, then we apply Lemma~\reft{cascade}, and
Lemma~\reft{A} to show that $l_{\Gamma}(h_g)\leq
(l_{\Gamma}(v)+C+c(2))+2C+2c(2)$. In any case, we get the
statement of the case~\refe{twob}.

Case \refe{gug} If $l_{\Gamma}(g_1)=l_{\Gamma}(g_3)$
(Figure~\ref{fig:fig582a}), then both $u$ and $v$ are conjugate to
$k\in H$ with $l_{\Gamma}(k)\leq C,$ and we obtain \refe{guga}.
Let $l_{\hat{\Gamma}}(\hat{g}_1)>l_{\hat{\Gamma}}(\hat{g}_3)$
(Figure~\ref{fig:fig582bc}, left). Assume that $\hat{p},$
$\hat{w_u}$ and $\hat{q}$ penetrate another coset $f_1H$ so that
$\hat{q}$ travels in $f_1H$ along $h_g.$ As both $\hat{w_u}$ and
$\hat{p}^{-1}$ are geodesics, the distances from their common
terminal point to a coset they both penetrate, are equal;
therefore, the case shown on Figure~\ref{fig:fig582bc}, left, is
the only possible one. In this case, the argument used to prove
Lemma~\reft{B}, implies that $l_{\Gamma}(h_g)\leq C.$ The other
cases (Figure~\ref{fig:fig582bc}, right, shows one of those) are
analogous to the case~\refe{twob}, so that $l_{\Gamma}(h_g)\leq
l_{\Gamma}(v)+3C+3c(2)$ which proves \refe{gugb}.
\begin{figure}[ht!]\centering
\includegraphics[width=.7\textwidth]{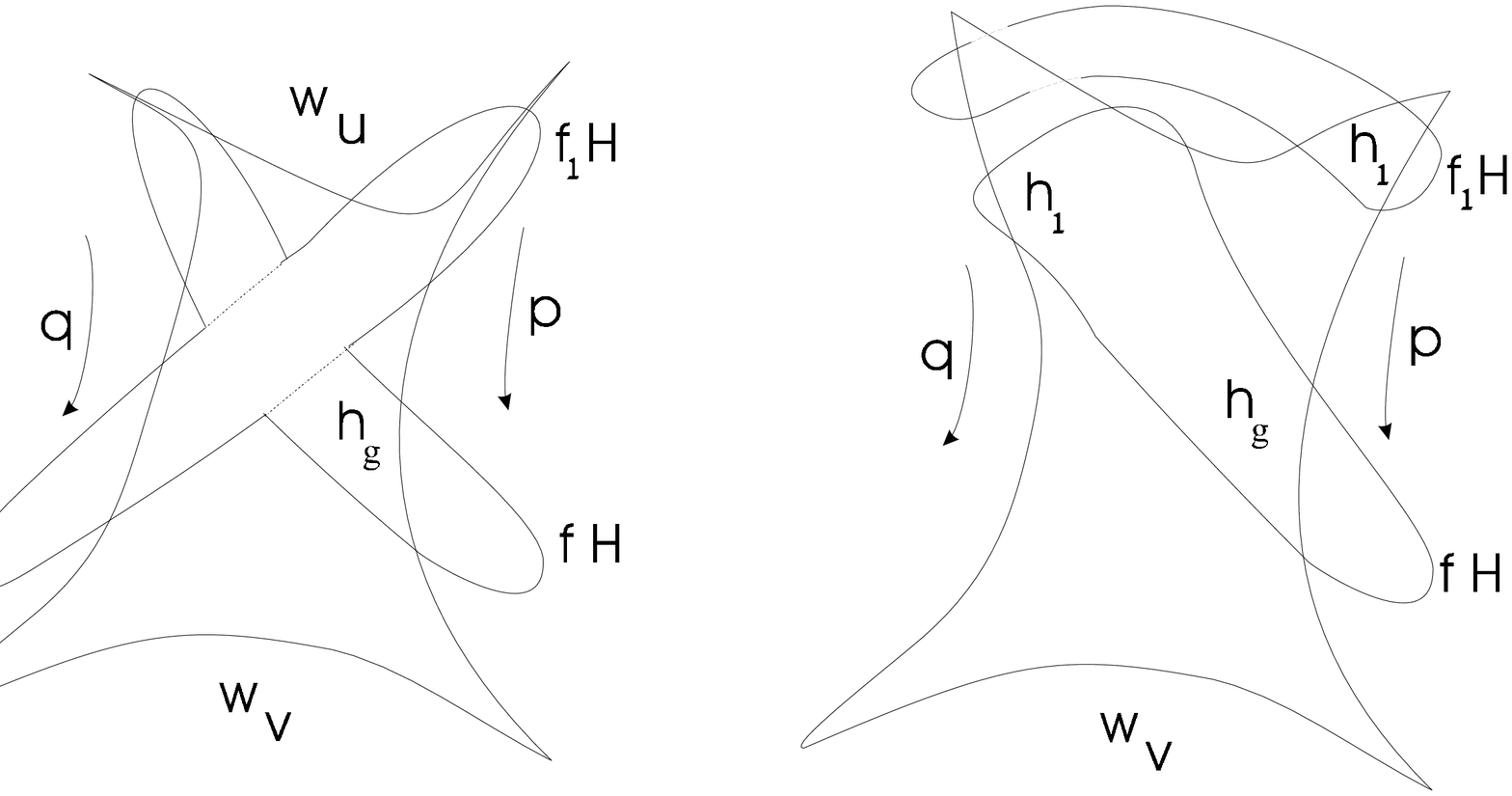}
\caption{Illustrated above is Lemma~\reft{mtl}, case~\refe{gug},
when $\hat{p},$ $\hat{w_u}$ and $\hat{q}$ penetrate two cosets of
$H$. \label{fig:fig582bc}}
\end{figure}
\begin{figure}[ht!]\centering
\includegraphics[width=.7\textwidth]{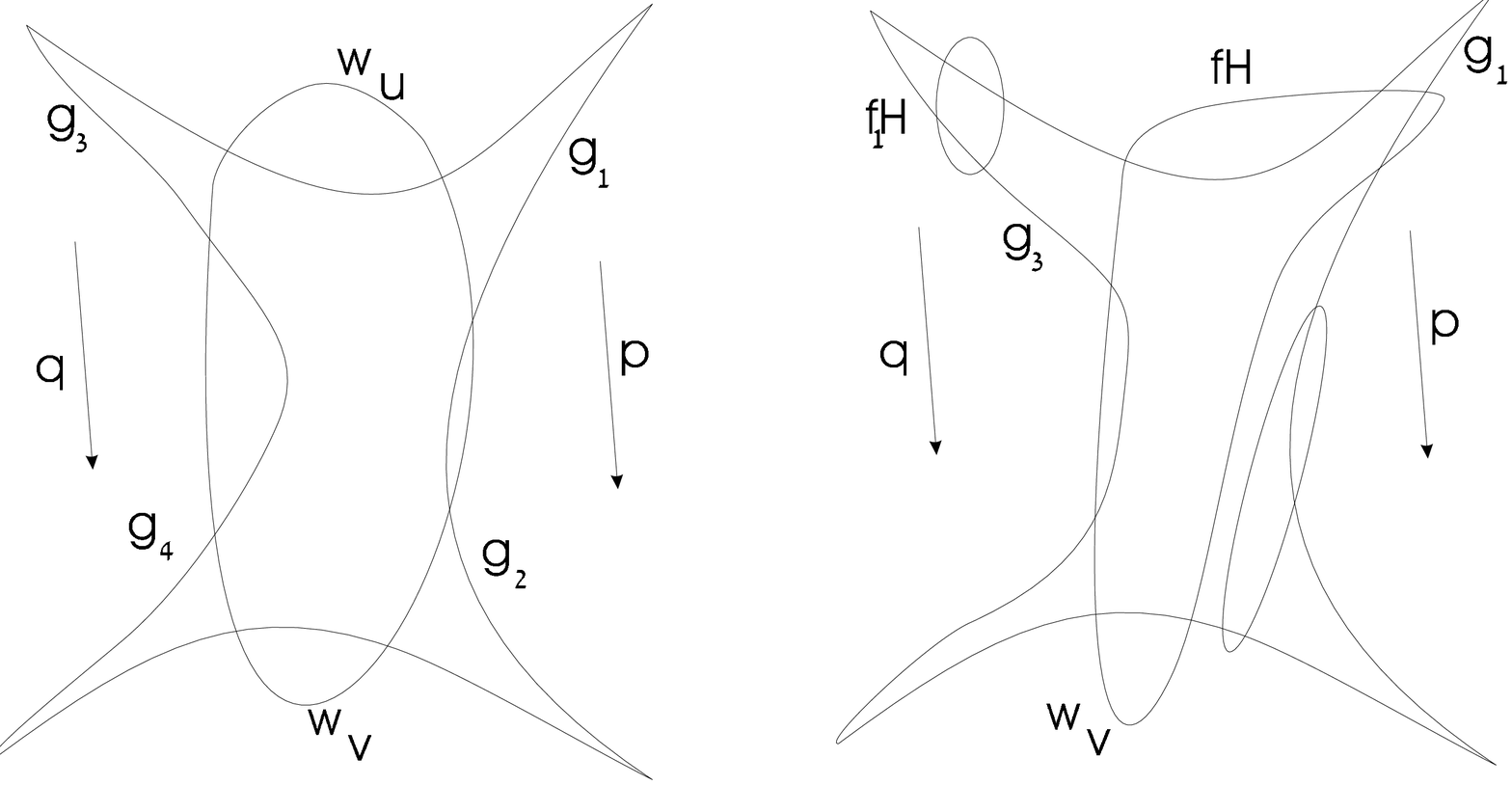}
\caption{Illustrated above is Lemma~\reft{mtl}, case~\refe{four},
when $\hat{p},$ $\hat{w}_u$, $\hat{q}$ and $\hat{w}_v$ penetrate a
coset of $H$. \label{fig:fig581}}
\end{figure}

Case \refe{four} follows easily from the above arguments. If
$l_{\hat{\Gamma}}(\hat{g}_1)=l_{\hat{\Gamma}}(\hat{g}_3)$
(Figure~\ref{fig:fig581}, left) then we have the
case~\refe{foura}. Indeed, it follows that $g_1=g_3$ and
$g_2=g_4.$ Let $k$ be a $\Gamma$-geodesic joining the points where
$p$ and $q$ leave $fH$; hence, both $u$ and $v$ are conjugate to
$k$ and $l_{\Gamma}(k)<l_{\Gamma}(v)+2C$, as claimed. In the other
cases (see, for instance, Figure~\ref{fig:fig582bc}, right, where
$l_{\hat{\Gamma}}(\hat{g}_1)<l_{\hat{\Gamma}}(\hat{g}_3)),$
Lemma~\reft{A} implies the statement of \refe{fourb}.
\end{proof}
\begin{corol}\label{t:comtl} If $ugvg^{-1}=1$ in $G,$
then either $u$ and $v$ are conjugate in $G$ to $k\in H$ with
$l_{\Gamma}(k)\leq Q+2C,$ or $g$ travels a $\Gamma$-distance
bounded by $2Q+ 6C$, in each coset it penetrates.
\end{corol}
\subsection{A global bound on the length of $g$}\label{s:global}
We have proved that if the relative length of $g$ is bounded, then
the $\Gamma$-length of $g$ can be bounded as well. \textit{A
priori}, we do not have any bound on the relative length of $g$.
It turns out that in order to bound globally the relative length
of $g$ (see the proof of Theorem~\reft{reLys} below), we need to
bound distances which $g$ travels in $H$-cosets. Let $D=c(8Q)$.
\begin{lem}\label{t:D} If $ugvg^{-1}=1$ in $G,$
then either
\begin{enumerate}
   \item \label{e:conjH} $u$ and $v$ are conjugate in $G$ to an element
   $k\in H$, or
   \item \label{e:nonconjH} The $\Gamma$-distance which $g$ travels in a coset
   it penetrates, is bounded by $l_H^g=2Q+10D$.
\end{enumerate}
\end{lem}
\begin{proof} Assume that the conditions of the case~\refe{conjH}
do not hold. Let $h_g$ be an $H$-subword of $g$ of the maximal
possible length.

First, assume that the path $h_g$ is \emph{not a subpath of an
$H$-floor of a cascade}. We claim that in this case the
$\Gamma$-length of $h_g$ satisfies the following inequality:
 \begin{equation}\label{e:nocascbound}
  l_{\Gamma}(h_g)\leq 2Q+6D.
 \end{equation}
Our proof of this latter claim splits according to the following
possibilities.
\begin{enumerate}
 \item \label{e:gshort} If $l_{\hat\Gamma}(\hat{g})< 3Q$, then
 $l_{\hat\Gamma}(\hat{w})< 8Q$, and by Corollary~\reft{comtl},
 we conclude that the inequality~\refe{nocascbound} holds.
 \item \label{e:glong} If $l_{\hat\Gamma}(\hat{g})\ge 3Q$, then
 by Lemma~\reft{lemmareL1}, $\lambda_1$ and $\lambda_2$
form a pair of $(2\hat{Q}+1)$-quasi-geodesics with common
endpoints. We distinguish the following two cases:
  \begin{enumerate}
  \item \label{e:nosamecoset} If $\hat{p}$ and $\hat{q}$ do not penetrate the same
    coset, then by Corollary~\reft{corolA}, $l_{\Gamma}(h_g)\leq
    2Q+c(2\hat{Q}+1)+2c(2).$ As $c(2\hat{Q}+1)<D$ and $c(2)<D$,
    it follows that $l_{\Gamma}(h_g)$ satisfies the inequality~\refe{nocascbound}.
  \item \label{e:nocasc} If $\hat{p}$ and $\hat{q}$ penetrate the same coset but
   $h_g$ is not a subpath of an $H$-floor of the cascade, then one can find a closed
   path $\hat{w}'$ which goes through $h_g$ and so that the subsegments of
   $\hat{p}$ and of $\hat{q}$ which belong to this closed path, do
   not penetrate the same coset. Therefore, the arguments used in
   the cases~\refe{gshort} and~\refe{nosamecoset} apply to $\hat{w}'$, so that the
   inequality~\refe{nocascbound} holds in this case also. Observe
   that if in $\hat{w}$ the cascade effect of length $n$ occurs, then $h_1$ and $h_{n+1}$
   occur also outside the cascade. Hence, the argument that we use in this case,
   applies to $h_1$ and $h_{n+1}$ as well. Thus, both $l_{\Gamma}(h_1)$ and $l_{\Gamma}(h_{n+1})$
   satisfy the inequality~\refe{nocascbound}.
  \end{enumerate}
\end{enumerate}

Now, assume that $h_g$ is a \emph{subpath of an $H$-floor of a
cascade}. By Lemma~\reft{cascade}, $l_{\Gamma}(h_g)\le
l_{\Gamma}(h_1)+2C_0+2c(2)$. Since $l_{\Gamma}(h_1)\leq 2Q+6D$
(see the the proof in the case~\refe{nocasc} above), and $C_0<D$,
we have that $l_{\Gamma}(h_g)\leq 2Q+10D$.
\end{proof}
The following theorem establishes explicitly the dichotomy
mentioned above: either $u$ and $v$ are conjugate to an element of
$H$, or the $\Gamma$-length of $g$ can be globally bounded.
\begin{thm}\label{t:reLys} Let $G$ be a group hyperbolic relative to
a subgroup $H,$ in the strong sense. If $ugvg^{-1}=1$ in $G$ and
the relative length of $g$ is positive and minimal possible, then
either
\begin{enumerate}
   \item \label{e:conH} $u$ and $v$ are conjugate in $G$ to an element
   $k\in H$ so that $g=g_1g_2$ where $u=g_1kg_1^{-1}$ and $k=g_2vg_2^{-1}$, or
   \item \label{e:relbounded} The $\Gamma$-length of $g$ is bounded
   in terms of the $\Gamma$-lengths of $u$ and $v$.
\end{enumerate}
\end{thm}
\begin{proof} We assume that the case~\refe{conH} does not occur.
  Let $\hat\gamma$ be a geodesic joining $\hat{p}(t)$ and
$\hat{q}(t)$, and let $m=l_0+2K$ (cf.\ \refe{gamma}) be the upper
bound for the $\hat{\Gamma}$-length of $\hat{\gamma}$ obtained in
the proof of Lemma \reft{reL2}. Assume that
$l_{\hat{\Gamma}}(\hat{g})>2\hat{Q}+2\delta+6m$ so that we can
consider values of $t$ satisfying the inequality $t_1+3m<t<t_2-3m$
$(t_1,t_2$ are as in Lemma \reft{reL2}). Denote by $\hat\gamma_1$
a geodesic joining $\hat{p}(t+3m)$ and $\hat{q}(t+3m),$ and denote
by $\hat\gamma_p$ (or $\hat\gamma_q)$ the segment of $\hat{p}$ (or
$\hat{q})$ between $\hat{p}(t)$ (or $\hat{q}(t)$) and
$\hat{p}(t+3m)$ (or $\hat{q}(t+3m)$). The closed path
$\hat{w}=\hat\gamma_p\hat{\gamma}_1\hat\gamma_q^{-1}\hat{\gamma}^{-1}$
has a relative length bounded by $8m$. Moreover, since
$l_{\hat{\Gamma}}(\hat{\gamma}),l_{\hat{\Gamma}}(\hat{\gamma_1})\leq
m$ and the distance between their initial (or terminal) points
equals $3m$, these two geodesics never penetrate the same coset.
Therefore, by Corollary~\reft{corolA}, the maximal distance that
$\hat{\gamma}$ (or $\hat{\gamma}_1$) can travel in a coset it
penetrates, is bounded above by $2l_H^g+3c(8m)$. Therefore, the
$\Gamma$-length of $\gamma$ is bounded in terms of the
$\Gamma$-lengths of $u$ and $v$, so that we can apply the argument
used to prove Theorem~\reft{hyperbolic}. This argument tells that
since $l_{\hat{\Gamma}}(\hat{g})$ is minimal possible, it is
bounded in terms of the $\Gamma$-lengths of $u$ and $v$. By
Lemma~\reft{D}, one obtains the claim.
\end{proof}
\subsection{The case when $u$ and $v$ are
conjugate to an element of $H$}\label{s:conjh} In general, there
is a finite sequence of elements $k_1,k_2,\dots,k_n$ of $H$
conjugate in $G$ to each other as follows:
$k_i=g_ik_{i+1}g_i^{-1}$, so that $u=g_uk_1g_u^{-1}$ and
$v=g_vk_ng_v^{-1}$ for some $g_i$, $g_u$ and $g_v$ in $G$. We are
able to find bounds for the $\Gamma$-length of $g_u$, of $g_v$ and
of those $g_i$ which are in $G\setminus H$, but if $g_i\in H$,
then its $\Gamma$-length cannot be bounded. Therefore, there is no
bound on the $\Gamma$-length of the element $g=g_ug_1\dots
g_{n-1}g_v$ conjugating $u$ and $v$. Note that in the geometric
picture, $\hat{p}$ and $\hat{q}$ penetrate the same cosets, at the
same moments of time, so that $\hat{w}$ is a finite sequence of
digons with two isosceles triangles at both ends of it. In the
case of a hyperbolic groups, this picture would mean that both $u$
and $v$ were conjugate to the trivial element, therefore, were
trivial elements themselves.

Our approach is as follows. Lemma~\reft{E} implies that the
$\Gamma$-length of $k_1$ and $k_n$ can be bounded in terms of the
length of $u$ and of $v$. By Corollary~\reft{hcofh},
$l_{\Gamma}(k_i)\leq c(2)$ for $1<i<n$. Moreover, the length of
$g_u$, of $g_v$ and of those conjugating elements $g_i$ which are
not in $H$, can be bounded in terms of the length of $u$ and of
$v$ as well. Therefore, it is enough to consider the finite set
$H_d$ of elements of $H$ whose length does not exceed $d=c(2)$.
Lemma~\reft{hhcon} below allows one to obtain the partition of
$H_d$ to conjugacy classes of $G$. Having obtained this partition,
we are able to establish whether or not $u$ and $v$ are conjugate
to each other, if each one of them is conjugate to an element of
$H$.
\begin{lem} \label{t:E} Let $u\in G$ be conjugate to $h\in H$.
Then either $u\in H$, and $u$ and $h$ are conjugate in $H$, or
there exist $k\in H$ and $g\in G$ so that $u=gkg^{-1}$ and the
following conditions hold:
\begin{enumerate}
  \item\label{e:E1} If $f\in G$ and $k_f\in H$ satisfy the equality $u=fk_ff^{-1}$,
  then $l_{\hat{\Gamma}}(\hat{g})\leq l_{\hat{\Gamma}}(\hat{f})$.
  \item\label{e:E2} $l_{\Gamma}(k)\leq C_0$.
  \item\label{e:E3} $l_{\Gamma}(g)$ is bounded in terms of the
  $\Gamma$-length of $u$.
\end{enumerate}
\end{lem}
\begin{proof} Assume that $u\notin H$. Corollary~\reft{membH} implies that
$g\in G\setminus H$, in particular the relative length of $g$ is
strictly positive. Therefore, the minimal possible relative length
is attained, and we get \refe{E1}. Fix an element $g$ that
satisfies the condition~\refe{E1}, and consider the closed path
$w=ugkg^{-1}$ and its projection $\hat{w}$ into $\hat{\Gamma}$. By
Lemma~\reft{lemmareL1}, either $l_{\hat{\Gamma}}(\hat{w})\leq
7\hat{Q}$, or $\hat{w}$ is the concatenation of two
$(2\hat{Q}+1)$-quasi-geodesics. The assertion~\refe{E2} follows
then either from Lemma~\reft{cboup}, or from the definition of the
BCP property. To obtain the assertion~\refe{E3}, apply
Theorem~\reft{reLys} and note that the assertion~\refe{E1} we have
just proven, implies that the case~\refe{conH} mentioned in the
statement of Theorem~\reft{reLys}, does not occur.
\end{proof}
\begin{corol}\label{t:cofh} Given $u\in G\setminus H,$ one can determine
effectively, whether or not there is $h\in H$ so that $u$ is
conjugate to $h.$
\end{corol}
\begin{proof} According to Lemma~\reft{E}, it is enough to determine
whether or not the word $ughg^{-1}$ is trivial for some $h$ and
$g$ whose $\Gamma$-length is bounded. There are only finitely many
possibilities to choose $h$ and $g$, and by \cite[Theorem
3.7]{Farb} (see Theorem~\reft{F37}), for each particular choice of
these elements, an answer can be found effectively.
\end{proof}
\begin{corol}\label{t:hcofh} If $u\in H$ and $h\in H$ are conjugate in $G$
but are not conjugate in $H$, then the assertion~\refe{E2} of
Lemma~\reft{E} becomes $l_{\Gamma}(k)\leq c(2)$.
\end{corol}
\begin{proof} Note that in this case the projection $\hat{w}$ of the
closed path $w=ugkg^{-1}$ is the concatenation of two
$2$-quasi-geodesics.
\end{proof}
\begin{lem} \label{t:hhcon} Given $h_u,h_v\in H,$ one can
determine effectively whether or not $h_u$ and $h_v$ are conjugate
in $G.$
\end{lem}
\begin{proof} First, we check whether or not $h_u$ and $h_v$ are conjugate
in $H.$ Assume that this is not the case. Let $d=c(2)$ be the
constant given by Definition~\reft{BCPp}. Consider the finite
subset $H_d=\{h\in H\mid l_{\Gamma}(h)\leq d\}$ of ``short"
elements of $H,$ and the partition of $H_d$ into conjugacy classes
$\mathcal{C}_G$ of $G:$ elements $h_1$ and $h_2$ of $H_d$ belong
to a $\mathcal{C}_G$-class if and only if there is $g\in G$ such
that $h_v=gh_ug^{-1}.$ We claim that this partition of $H_d$ can
be obtained in a finite time. Indeed, as the conjugacy problem in
$H$ is solvable, we can find a partition $\mathcal{C}_H$ of $H_d$
into conjugacy classes of $H$ in a finite time. Furthermore, we
define bounded $\mathcal{C}_G$-classes as follows: elements
$\tilde{k}$ and $\tilde{h}$ of $H_d$ belong to a bounded
$\mathcal{C}_G$-class if and only if either $\tilde{k}$ and
$\tilde{h}$ belong to a $\mathcal{C}_H$-class, or there is a
finite sequence $\tilde{k}=k_1,k_2,\dots,k_n=\tilde{h}$ of
elements of $H$ with $l_{\Gamma}(k_j)\leq c(2)$, so that for every
$i=1,2,\dots,n-1$ there is $g_i\in G$ with bounded length
$l_{\Gamma}(g_i)$ such that $h_i=g_ih_{i+1}g_i^{-1}.$
Corollary~\reft{hcofh} implies that two elements of $H_d$ belong
to a $\mathcal{C}_G$-class if and only if they belong to a bounded
$\mathcal{C}_G$-class. This observation gives rise to the
following algorithm. Pick a $\mathcal{C}_H$-class
$H_d^{(1)}=\{h_1^{(1)},\dots,h_{m_1}^{(1)}\}$. For each
$h_i^{(1)}\in H_d^{(1)}$, find all those pairs of elements $k\in
H$ and $g\in G\setminus H$ which satisfy the conditions of
Lemma~\reft{E} and Corollary~\reft{hcofh}. Since each $k$ is a
``short" element of $H$, it belongs to a $\mathcal{C}_H$-class
$H_d^{(j)}$. We add all these classes $H_d^{(j)}$ to $H_d^{(1)}$
so as to obtain a bounded $\mathcal{C}_G$-class, and declare all
these added elements as new members in the $\mathcal{C}_G$-class.
Having collected all $k$ and their $\mathcal{C}_H$-classes, we
repeat the above procedure for each new member in the
$\mathcal{C}_G$-class of $h_u.$ Again, added elements are declared
to be new members, and we proceed with them in the same manner,
until there are no new members anymore. Then we pick a
$\mathcal{C}_H$-class, which is not a subset of the bounded
$\mathcal{C}_G$-class of $h_u$ we have just obtained, and repeat
the same procedure.
 The algorithm stops when the (finite) set of $\mathcal{C}_H$-classes is
exhausted.

Having obtained the partition of $H_d$ into
$\mathcal{C}_G$-classes, we check whether or not there are
$k_u,k_v\in H_d$ that belong to same $\mathcal{C}_G$-class and
such that $k_u$ and $h_u$ as well as $k_v$ and $h_v,$ are
conjugate in $H$ by elements of bounded length (see the
assertion~\refe{E3} of Lemma~\reft{E}). The elements $h_u$ and
$h_v$ are conjugate in $G$ if and only if there are $k_u$ and
$k_v$ as above.
\end{proof}
From Corollary~\reft{cofh} and Lemma~\reft{hhcon} we obtain the
following corollary.
\begin{corol}\label{t:uhcon} Given $u\in G\setminus H$ and $h\in H,$
one can determine effectively, whether or not $u$ and $h$ are
conjugate in $G.$
\end{corol}
\subsection{Proof of Theorem~\reft{conj}}
\begin{proof} By Corollary~\reft{membH}, we can determine whether or
not $u$ and $v$ belong to $H$. If both $u$ and $v$ are in $H,$
then the assertion follows from Lemma \reft{hhcon}. If for
instance $u\in G\setminus H$ while $v\in H,$ then the assertion
follows from Corollary~\reft{uhcon}. Now, assume that neither $u$
nor $v$ is in $H.$ According to Corollary \reft{cofh}, we can
answer effectively the following two questions:
\begin{enumerate}
 \item Is there $k_u\in H$ so that $k_u$ and $u$ are conjugate in
 $G?$
 \item Is there $k_v\in H$ so that $k_v$ and $v$ are conjugate in
 $G?$
\end{enumerate}
If the answers are different, then $u$ and $v$ are not conjugate
in $G.$ If both answers are positive, then we apply Lemma
\reft{hhcon} to $k_u$ and $k_v;$ $u$ and $v$ are conjugate in $G$
if and only if $k_u$ and $k_v$ are conjugate in $G.$ If both
answers are negative, then by Theorem~\reft{reLys}, $u$ and $v$
are conjugate if and only if there is a conjugating element of a
bounded $\Gamma$-length. Since balls of bounded radii in the
Cayley graph $\Gamma$ of $G$ are compact, this latter condition
can be checked effectively.
\end{proof}
\subsection{Group with several parabolic subgroups}
The definition of a relatively hyperbolic group can be extended to
the case of several subgroups \cite[Section 5]{Farb}. Let $G$ be a
group, and let $\{H_1,\dots,H_r\}$ be a finite set of finitely
generated subgroups of $G.$ In the Cayley graph of $G,$ for every
$i=1,2,\dots,r,$ add a vertex $v(gH_i)$ for each left coset of $H_i$
in $G,$ and connect this new vertex (by an edge with length $\frac
12,$ as before) with each element of this left coset. This new graph
$\hat{\Gamma}$ is called the {\it coned-off graph} of $G$ with
respect to $\{H_1,\dots,H_r\}.$ The group $G$ is {\it weakly
hyperbolic relative to} $\{H_1,\dots,H_r\},$ if $\hat{\Gamma}$ is a
hyperbolic metric space. The definition of the BCP property can be
extended in an obvious way to this case. If the subgroups
$H_1,\dots,H_r$ are torsion-free, then the BCP property implies that
these subgroups are pairwise conjugacy separated. This means that if
$gH_ig^{-1}\cap H_j\ne\emptyset$ for some $g\in G,$ $1\le i,j\le r,$
then necessarily $i=j$ and $g\in H_i.$ The group $G$ is
\emph{strongly hyperbolic relative to the family of subgroups}
$\{H_1,\dots,H_r\},$ if $G$ is weakly hyperbolic relative to
$\{H_1,\dots,H_r\},$ and the pair $(G,\{H_1,\dots,H_r\})$ has the
BCP property.

Our arguments can easily be extended to prove the following
generalization of Theorem \reft{conj}.
\begin{thm}\label{t:confs} Let $G$ be a group strongly hyperbolic relative
to a finite set of subgroups $\{H_1,\dots,H_r\}$. If the conjugacy
problem is solvable in $H_i$ for all $i=1,2,\dots,r,$ then it is
solvable in $G.$
\end{thm}

\section{Fundamental groups of negatively curved
manifolds}\label{s:manifolds} In this section, we prove Theorem
\reft{mnfld}. Let $G$ be the fundamental group of a negatively
curved non-compact manifold of finite volume with a single cusp,
and let $H$ denote the cusp subgroup of $G.$ Then $G$ is
hyperbolic relative to $H$ in the strong sense (Farb~\cite{Farb}
gave a direct proof of this assertion). Since $H$ is a nilpotent
group, the conjugacy problem for $H$ is solvable \cite{M}.
Therefore, in order to prove Theorem~\reft{mnfld}, it remains to
show that the constants $c(P)$ in Definition \reft{BCPp} can be
bounded effectively.

We follow \cite{Farb} and \cite{Reb}. $\tilde{M}$ denotes a
Hadamard manifold; we are most interested in the case when
$\tilde{M}$ is the universal cover of a complete, finite-volume
negatively curved Riemannian manifold $M$ with pinched negative
curvature $-b^2\le K(M)\le -a^2<0.$ Our calculation is based on
the geometry of horospheres in $\tilde{M}.$ Let $x\in \tilde{M},$
$z$ be a point at infinity, and let $\gamma$ be the geodesic ray
from $x$ to $z.$ {\it Horospheres} are the level surfaces of the
Busemann function $F=\lim_{t\rightarrow\infty} F_t,$ where $F_t$
is defined by $F_t(p)=d_{\tilde{M}}(p,\gamma(t))-t.$ Let $S$ be a
horosphere, we denote by $d_S$ the induced path metric on $S;$
that is, $d_S(x,y)$ is the infimum of the length of all paths in
$S$ from $x$ to $y.$

\begin{prop}\label{t:F42}{\rm\cite[Proposition 4.2]{Farb}}\qua
If $\gamma$ is a geodesic tangent to $S,$ and $p$ and $q$ are
projections of $\gamma(\pm \infty)$ onto $S,$ then
$$\frac 2b\le d_S(p,q)\le \frac 2a$$
\end{prop}

\begin{prop}[{\cite[Corollary 5.3]{Reb}, cf.\
\cite[Proposition 4.3]{Farb}}]\label{t:F43}  
Let $S$ and $S'$ be nonintersecting
horospheres based at distinct points of $\partial H.$ Then the
$S$-diameter of the projection $\pi_S(S')$ is at most $\frac
4a+2\delta,$ where $\delta$ is the Gromov hyperbolicity constant
for $\tilde{M}.$
\end{prop}

\begin{defn}\label{t:vsh}\cite{Farb}\qua
Let $\gamma$ be a geodesic in $\tilde{M}$ not intersecting a
horosphere $S$. Given $s\in S$, we say that $\gamma$ {\it can be
seen from} $s$ if $\overline{s\gamma(t)}\cap S=\{s\}$ for some
$t.$ Let $T_{\gamma}$ be the set of points $s\in S$ that $\gamma$
can be seen from. The {\it visual size} $V_S$ of the horosphere
$S$ is defined to be the supremum of the diameter of $T_{\gamma}$
in the metric $d_S,$ where the supremum is taken over all
geodesics $\gamma$ not intersecting $S.$
\end{defn}
\begin{prop}[{\cite[Lemma 5.4]{Reb}, cf.\
\cite[Lemma 4.4]{Farb}}]\label{t:F44}  
Horospheres in a pinched Hadamard manifold
have (uniformly) bounded visual size.
\end{prop}
The proof shows that the visual size of $S$ is bounded by $\frac
2a+C$ where according to \cite[Lemma 4.10]{Reb}, $C=2\delta+\log
16.$ Therefore, the visual size $V_S$ of $S$ satisfies the
following inequality:
\begin{equation}\label{e:vsh}
V_S\le\frac 2a+2\delta+\log 16.
\end{equation}
Let $G$ be the fundamental group of $M$ so that $M=\tilde{M}/G.$
We can choose a $G$-invariant set of horoballs so that there is a
uniform lower bound on the distance between horoballs and the
action of $G$ on the horoballs has finitely many orbits. Having
deleted the interiors of all of these horoballs, we obtain a space
$X$ on which $G$ acts cocompactly by isometries $(X$ is equipped
with the path metric). Choose a base point $x\in X,$ the map
$g\mapsto g\cdot x$ gives a quasi-isometry $\h{\psi}{\Gamma}X$ of
the Cayley graph of $G$ with $X;$ for each coset $gH,$ all of the
elements of $gH$ are mapped to the same horosphere. The {\it
electric space} $\hat{X}$ is the quotient of $X$ obtained by
identifying points which lie in the same horospherical boundary
component of $X.$ The quotient $\hat{X}$ has a path pseudo-metric
$d_{\hat{X}}$ induced from the path metric $d_X;$ the
pseudo-metric $d_{\hat{X}}$ can be thought of as a pseudo-metric
on $X,$ where the distance between two points is the length of the
shortest path between them, but path-length along a horosphere
$S\subset\partial X$ is measured as zero length. Locally
$d_{\hat{X}}$ agrees with $d_{\tilde{M}}$ on the interior of $X.$
\begin{prop}\label{t:F46}{\rm\cite[Proposition 4.6]{Farb}}\qua
The electric space $\hat{X}$ is a $\delta'$-hyperbolic
pseudometric space for some $\delta'>0.$
\end{prop}
Given a path $\gamma$ in $\hat{X},$ the {\it electric length}
$l_{\hat{X}}(\gamma)$ is the sum of the $X$-length of subpaths of
$\gamma$ lying outside every horosphere. An {\it electric
geodesic} between $x,y\in\hat{X}$ is a path $\gamma$ in $\hat{X}$
from $x$ to $y$ such that $l_{\hat{X}}(\gamma)$ is minimal. An
{\it electric $P$-quasi-geodesic} is a $P$-quasi-geodesic in
$\hat{X}.$

\begin{lem}[{\cite[Lemma 5.6]{Reb}, cf.\
\cite[Lemma 4.5]{Farb}}] \label{t:F45} Given $P>0,$ there exist constants
$K=K(P),$ $L=L(P)>0$ such that for any electric $P$-quasi-geodesic
$\beta$ from $x$ to $y,$ if $\gamma$ is the $\tilde{M}$ geodesic
from $x$ to $y,$ then $\beta$ stays completely inside
$Nbhd_{\hat{X}}(\gamma,K+L/2).$
\end{lem}
According to the proof, $K$ can be chosen so that
\begin{equation}\label{e:K45}
K\ge\frac 1a\log(2P(V_S+1)).
\end{equation}
Then one can set
\begin{equation}\label{e:L45}
L=4PK(2+V_S)+8P\delta,
\end{equation}
where $V_S$ is as in \refe{vsh}.
\begin{lem}[{\cite[Lemma 5.7]{Reb}, cf.\
\cite[Lemma 4.7]{Farb}}] \label{t:F47} Let $\beta$ be an electric
$P$-quasi-geodesic so that $\beta\cap S=\emptyset.$ Then there
exists a constant $D=D(P)$ so that $\pi_S(\beta)$ has $S$-length
at most $Dl_{\hat{X}}(\beta).$
\end{lem}
The proof shows that
\begin{equation}\label{e:F47}
D=1+V_S\le 1+\frac 2a+2\delta+\log 16.
\end{equation}
\begin{lem}[{\cite[Lemma 5.9]{Reb}, cf.\
\cite[Lemma 4.8 and Lemma 4.9]{Farb}}]\label{t:F48}  
Let $\alpha$ and $\beta$ be
electric $P$-quasi-geodesics from $x$ to $y$ in $\hat{X}.$ Then
there exists a constant $E$ such that the following conditions
hold:
\begin{enumerate}
 \item Suppose $\alpha$ first greets $S$ at $\alpha(s_0)$ and
$\beta$ first greets $S$ at $\beta(t_0).$ Suppose that $\alpha$ and
$\beta$ leave $S$ at $\alpha(s_1)$ and $\beta(t_1).$ Then
\[
d_S(\alpha(s_0),\beta(t_0))<E\quad \text{and}\quad
d_S(\alpha(s_1),\beta(t_1))<E.
\]
 \item Suppose $\alpha$ greets $S$
at $\alpha(s_0)$ and leaves $S$ at $\alpha(s_1).$ Suppose that
$\beta$ doesn't greet $S.$ Then $d_S(\alpha(s_0),\alpha(s_1))<E.$
\end{enumerate}
\end{lem}
By the proof, $E=3\max\{\delta,D\}P(K+L),$ where $D,$ $K$ and $L$
are given by \refe{F47}, \refe{K45} and \refe{L45}, respectively.
Altogether, we have
\begin{equation}\label{e:bhp}
E\le 3P(1+V_S)(K+L),
\end{equation}
where $V_S$ is as in \refe{vsh}. Therefore, the upper bound for
the constant $E$ can be computed effectively. Let $\lambda$ denote
the quasi-isometry constant of the map $\h{\psi}{\Gamma}X.$ Then
the constants $c(P)$ of Definition \reft{BCPp} can be bounded as
follows: $c(P)\le \lambda E.$
\rk{Acknowledgements}
The first version of this paper was written when I was a
postdoctoral fellow at the Technion, Haifa. I am thankful to Arye
Juh\'{a}sz for stimulating discussions, and to Benson Farb for
many helpful electronic communications. Also, I want to thank
Michah Sageev for pointing me out the example of a weakly
relatively hyperbolic group with unsolvable conjugacy problem.
Finally, I am deeply grateful to the referee for extremely useful
remarks and helpful suggestions.

\Addresses\recd


\begin{thebibliography}

\bibitem{Al} {\bf E Alibegovic}, {\em A combination theorem for relatively hyperbolic groups}, Bull. London Math. Soc. (to appear)

\bibitem{ABC} {\bf J\,M Alonso}, {\bf T Brady}, {\bf D Cooper et al},
{\em Notes on word hyperbolic groups}, from: ``Group theory from a
geometrical viewpoint'', (E Ghys, A Haefliger, A
Verjovsky, editors) ICTP, Trieste (1990) 3--63
\MR{1170362}

\bibitem{Ba} {\bf B Baumslag}, {\em Residually free groups},
Proc. London Math. Soc.  17 (1967) 402--418 \MR{0215903}

\bibitem{Bow} {\bf B\,H Bowditch}, {\em Relatively hyperbolic groups},
preprint, University of Southampton (1998)

\bibitem{Bow0} {\bf B\,H Bowditch}, {\em Connectedness properties of limit sets},
Trans. Amer. Math. Soc. 351 (1999) 3673--3686 \MR{1624089}

\bibitem{Bow1} {\bf B\,H Bowditch}, {\em Boundaries of geometrically finite
groups}, Math. Z. 230 (1999) 509--527 \MR{1680044}

\bibitem{Bow2} {\bf B\,H Bowditch}, {\em Peripheral splittings of groups},
Trans. Amer. Math. Soc. 353 (2001) 4057--4082 \MR{1837220}

\bibitem{Bu} {\bf I Bumagin}, {\em On definitions of relatively
hyperbolic groups}, Proc. Amer. Math.
Soc. (to appear)

\bibitem{Can} {\bf J Cannon}, {\em The combinatorial structure of
cocompact discrete hyperbolic groups}, Geometriae Dedicata 16
(1984) 123--148
\MR{0758901}

\bibitem{Mi} {\bf D Collins}, {\bf C Miller}, {\em The conjugacy problem
and subgroups of finite index}, Proc. London Math. Soc. 34 (1977)
535--556
\MR{0435227}

\bibitem{Dacl} {\bf F Dahmani}, {\em Classifying space and boundary for
relatively hyperbolic groups}, Proc. London Math. Soc.(3) 86 (2003)
666--684 \MR{1974394}

\bibitem{Daco} {\bf F Dahmani}, {\em Combination of convergence groups}, 
\gtref7{2003}{27}{933}{963} \MR{2026551}

\bibitem{Dehn} {\bf M Dehn}, {\em \"{U}ber unendliche diskontinuierliche
Gruppen}, Math. Ann. 71 (1912) 116--144

\bibitem{Farb} {\bf B Farb}, {\em Relatively hyperbolic groups},
Geom. Func. Anal. 8 (1998) 810--840
\MR{1650094}

\bibitem{J} {\bf A Juh\'{a}sz}, {\em Extension of group presentations and
relative small cancellation theory. I}, Internat. J.
Algebra Comput. 10 (2000) 375--398
\MR{1759186}

\bibitem{GH} {\bf E Ghys}, {\bf P de\,la\,Harpe, Editors},
{\em Sur les groupes hyperboliques d'apr\`{e}s Mikhael Gromov},
Progress in Mathematics 83, Birkhauser (1990)
\MR{1086648}

\bibitem{Go} {\bf B Goldfarb}, {\em Novikov conjectures and relative
hyperbolicity}, Math. Scand. 85 (1999) 169--183 \MR{1724233}

\bibitem{Gr} {\bf M Gromov}, {\em Hyperbolic groups},
from: ``Essays in group theory'', (S~Gersten, editor)
Math. Sci. Res. Inst. Publ. 8, 
Springer (1987)  75--263
\MR{0919829}

\bibitem{Hr} {\bf G\,C Hruska}, {\em Nonpositively curved 2-complexes with isolated flats},
\gtref8{2004}{5}{205}{275}
\MR{2033482}

\bibitem{KS} {\bf I Kapovich}, {\bf P Schupp},
{\em Relative hyperbolicity and Artin groups}, Geom.
Dedicata (to appear)

\bibitem{KM} {\bf O\,G Kharlampovich}, {\bf A\,G Myasnikov},
{\em Irreducible affine varieties over a free group II}, J. of
Algebra 200 (1998) 517--570
\MR{1610664}

\bibitem{KMRS} {\bf O\,G Kharlampovich}, {\bf A\,G Myasnikov}, {\bf
V\,N Remeslennikov}, {\bf D\,E Serbin}, {\em Subgroups of fully
residually free groups: algorithmic problems}, Contemp. Math. 360
(2004) 63--101

\bibitem{Lys} {\bf I Lys\"{e}nok}, {\em Some algorithmic properties of hyperbolic
groups} (Russian), Izv. Akad. Nauk SSSR Ser. Mat. 53 (1989)
814--832; translation in Math. USSR-Izv. 35 (1990) 145--163
\MR{1018749}

\bibitem{MaM} {\bf H Masur}, {\bf Y Minsky}, {\em Geometry of the
complex of curves I: Hyperbolicity}, Invent. Math.  138 (1999) 103-149
\MR{1714338}

\bibitem{M} {\bf A Mostowski}, {\em On the decidability of some problems in
special classes of groups}, Fund. Math. 59 (1966) 123--135
\MR{0224693}

\bibitem{MRS} {\bf A\,G Myasnikov}, {\bf V\,N Remeslennikov}, {\bf
D\,E Serbin}, {\em Regular free length functions on Lyndon's free
$\mathbb{Z}[t]$-group $F^{\mathbb{Z}[t]}$}, Contemp.
Math. (to appear)

\bibitem{Osin} {\bf D Osin}, {\em Relatively hyperbolic groups: Intrinsic geometry, algebraic
properties, and algorithmic problems}, Mem. Amer. Math. Soc. (to appear)

\bibitem{Osin2} {\bf D Osin}, {\em Weak hyperbolicity and free
constructions}, Contemp. Math. (to appear)

\bibitem{Reb} {\bf D Rebbechi}, {\em Algorithmic Properties of Relatively Hyperbolic Groups},
Ph.D. thesis, Rutgers (2000)

\bibitem{Sd} {\bf Z Sela}, {\em Diophantine geometry over groups I},
Publ. Math. Inst. Hautes \'Etudes Sci. 93 (2001) 31-105 
\MR{1863735}

\bibitem{Sz} {\bf A Szczepanski}, {\em Relatively hyperbolic groups},
Michigan Math. J. 45 (1998) 611--618
\MR{1653287}

\bibitem{Sz2} {\bf A Szczepanski}, {\em Examples of relatively
hyperbolic groups}, Geom. Dedicata 93 (2002) 139--142
\MR{1934694}

\bibitem{Y} {\bf A Yaman}, {\em A topological characterization of
relatively hyperbolic groups}, J. Reine Angew. Math. 566 (2004)
41--89
\MR{1934694}

\end{thebibliography}
\end{document}